\documentclass[11pt, leqno, twoside]{article}

\usepackage{amsfonts}
\usepackage{mathrsfs}
\usepackage{amssymb}
\usepackage{amsmath}

\allowdisplaybreaks

\def\rr{{\mathbb R}}
\def\rn{{{\rr}^n}}

\def\cc{{\mathbb C}}
\def\nn{{\mathbb N}}
\def\zz{{\mathbb Z}}

\def\cb{{\mathcal B}}

\def\cj{{\mathcal J}}

\def\cm{{\mathcal M}}

\def\cx{{\mathcal X}}
\def\cy{{\mathcal Y}}

\def\noz{\nonumber}
\def\fz{\infty}

\def\supp{{\mathop\mathrm{\,supp\,}}}

\def\loc{{\mathop\mathrm{\,loc\,}}}

\def\BMO{{\mathop\mathrm{ BMO }}}

\def\blo{{\mathrm {\,BLO}}}

\def\essinf{{\mathop{\mathrm{\,essinf\,}}}}

\def\lz{\lambda}

\def\dz{\delta}

\def\ez{\epsilon}

\def\bz{\beta}
\def\gz{{\gamma}}

\def\kz{{\kappa}}
\def\tz{\theta}
\def\sz{\sigma}
\def\wz{\widetilde}

\def\ls{\lesssim}

\def\pz{{\prime}}

\def\gfz{\genfrac{}{}{0pt}{}}

\def\hs{\hspace{0.3cm}}

\def\dint{\displaystyle\int}
\def\dfrac{\displaystyle\frac}

\def\r{\right}
\def\lf{\left}

\newtheorem{thm}{Theorem}[section]
\newtheorem{lem}{Lemma}[section]

\newtheorem{cor}{Corollary}[section]

\newtheorem{pf}{Proof.}

\numberwithin{equation}{section}

\begin{document}

\arraycolsep=1pt

\title{{\vspace{-5cm}\small\hfill\bf Israel J. Math. (to appear)}\\
\vspace{4cm}\bf\Large
Characterizations of BMO Associated with Gauss Measures
via Commutators of Local Fractional
Integrals
\footnotetext{\hspace{-0.22cm}2000 {\it Mathematics Subject
Classification}. Primary 47B47; Secondary 47H50, 42B30.
\endgraf {\it Key words and phrases}.
Gauss measure, $\BMO(\gz)$, fractional integral operator,
fractional maximal operator,
commutator.
\endgraf The second author is supported by the National
Natural Science Foundation (Grant No. 10871025) of China.
\endgraf $^\ast$Corresponding author.}}
\author{Liguang Liu and Dachun Yang $^\ast$}
\date{ }
\maketitle

\begin{center}
\begin{minipage}{14cm}
\small\noindent{\bf Abstract } Let
$d\gz(x)\equiv\pi^{-n/2}e^{-|x|^2}dx$ for all $x\in\rn$
be the Gauss measure on $\rn$.
In this paper, the authors establish the characterizations of the
space $\BMO(\gz)$ of Mauceri and Meda via commutators of
either local fractional integral operators or local
fractional maximal operators. To this end, the authors first
prove that such a local fractional integral operator of order $\bz$
is bounded from $L^p(\gz)$ to $L^{p/(1-p\bz)}(\gz)$,
or from the Hardy space $H^1(\gz)$ of Mauceri and Meda to $L^{1/(1-\bz)}(\gz)$
or from  $L^{1/\bz}(\gz)$ to $\BMO(\gz)$,
where $\bz\in(0, 1)$ and $p\in(1, 1/\bz)$.
\end{minipage}
\end{center}

\section{Introduction\label{s1}}

\hskip\parindent The space $\BMO(\rn)$
of functions with bounded mean oscillation was first introduced by
John and Nirenberg \cite{jn} and it plays an important
role in harmonic analysis and partial differential equations;
see, for example, \cite{tor, g, st93}. One of the remarkable
characterizations of the space $\BMO(\rn)$ is given in terms of
commutators of certain operators. In particular, when $T$ is
a singular integral with standard kernel, Coifman, Rochberg
and Weiss \cite{crw} proved that $b\in \BMO(\rn)$ is sufficient for
$[b, T](f)\equiv bT(f)-T(bf)$  to be bounded on $L^p(\rn)$ with
$p\in(1, \fz)$ and also established a partial converse.
The full converse of this result was
obtained by Janson \cite{j78}. Moreover, assuming that $I_\bz$ is a fractional
integral operator of order $\bz$ with $\bz\in(0, n)$
(see, for example, \cite[p.\,116]{st70}), Chanillo \cite{ch82} proved that $[b, I_\bz]$
is bounded from $L^p(\rn)$ to $L^q(\rn)$ if and only if
$b\in\BMO(\rn)$, where $1<p<q<\fz$ and $1/q=1/p-\bz/n$.

The main purpose of this paper is to generalize the above result of
Chanillo \cite{ch82} to the setting of the Gauss measure metric
space $(\rn, |\cdot|, d\gz)$, where $|\cdot|$ denotes the Euclidean
norm and $d\gz(x)\equiv\pi^{-n/2}e^{-|x|^2}dx$ for all $x\in\rn$ the
Gauss measure. Such an underlying space $(\rn, |\cdot|, d\gz)$
naturally appears in the analysis associated with the
Ornstein-Uhlenbeck operator; see \cite{pi88,  g94, s97, gmst99, p01,
mm} and the references therein. However, $(\rn, |\cdot|, d\gz)$ is
not a space of homogeneous type in the sense of Coifman and Weiss
\cite{cw1,cw2}. Recently, Mauceri and Meda \cite{mm} developed a
theory of singular integrals on $(\rn, |\cdot|, d\gz)$ which plays
for the Ornstein-Uhlenbeck operator the same role as that the theory
of classical Calder\'on-Zygmund operators plays for the Laplacian on
classical Euclidean spaces. The approach used in \cite{mm} requires
the introduction of certain Hardy space $H^1(\gz)$ and its dual
space $\BMO(\gz)$ related to a certain class $\cb_a$ with $a\in(0,
\fz)$ of admissible balls.

In this paper, we characterize
the space $\BMO(\gz)$ of Mauceri and Meda \cite{mm}
via commutators of either local fractional integral operators
or local fractional maximal operators.
To this end, we first establish the boundedness of
such a local fractional integral operator on Lebesgue spaces
or the corresponding Hardy space and its dual space.
A main difficulty to obtain these results exists in
the non-doubling property of the Gauss measure.

To state our results, we first recall some notation and notions;
see, for example, \cite{mm}. Let $m(x)\equiv\min\{1, 1/|x|\}$
for all $x\in\rn$. We denote by $c_B$ and $r_B$, respectively, the center
and radius of any ball $B$. For any $\kz>0$,
denote by $\kz B$ the ball with center $c_B$ and radius
$\kz r_B$. Let $a\in(0,\fz)$. The admissible class $\cb_a$ of
balls is defined to be the set of all balls $B\subset\rn$ such that
$ r_B\le am(c_B)$. For any $a\in(0, \fz)$ and $x\in\rn$, denote by
$\cb_a(x)$ the collection of all balls $B\in\cb_a$ containing $x$.

Mauceri and Meda \cite[p.\,281]{mm} introduced the following
$\BMO(\gz)$ space. Precisely, a function $f\in L^1(\gz)$ is said to
belong to the space $\BMO(\gz)$ if
$$\|f\|_\ast\equiv\sup_{B\in\cb_1}
 \dfrac1{\gz(B)}\dint_B|f(x)-f_B|\,d\gz(x)<\fz,$$
where and in what follows, $\gz(B)$ denotes the Gauss measure
of $B$ and
$$f_B\equiv \dfrac1{\gz(B)}\dint_Bf(y)\,d\gz(y).$$
Moreover, the $\BMO(\gz)$ norm of $f$ is defined by
\begin{equation}\label{1.1}
\|f\|_{\BMO(\gz)}\equiv\|f\|_\ast+\|f\|_{L^1(\gz)}.
\end{equation}

Mauceri and Meda \cite{mm} also introduced the
atomic Hardy space $H^1(\gz)$, which is the predual space of $\BMO(\gz)$;
see \cite[Theorem~5.2]{mm}. Precisely,
assume $r\in(1, \fz)$. A $(1, r)$ atom is either the constant
function $1$, or a function $a\in L^1(\gz)$ supported in a ball $B\in\cb_1$
with the properties $\|a\|_{L^r(\gz)}\le[\gz(B)]^{1/r-1}$
and $\int_B a(x)\,d\gz(x)=0$.
The Hardy space $H^{1,\,r}(\gz)$ is the space of all functions $g\in L^1(\gz)$
that admits a decomposition of the form
\begin{equation}\label{1.2}
g=\sum_{k=1}^\fz\lz_ka_k,\end{equation}
where $\{a_k\}_{k=1}^\fz$ are $(1, r)$ atoms and
$\{\lz_k\}_{k=1}^\fz\subset\cc$ satisfying that $\sum_{k=1}^\fz|\lz_k|<\fz$.
The norm $\|g\|_{H^{1,\,r}(\gz)}$ of $g$ is defined to be
the infimum of $\sum_{k=1}^\fz|\lz_k|$ over all decompositions of $g$ as
in \eqref{1.2}. It was pointed out in \cite[p.\,297]{mm}
that the Hardy spaces $H^{1,\,r}(\gz)$
for all $r\in(1, \fz)$ coincide with equivalent norms,
which will be simply denoted by $H^1(\gz)$.

Motivated by geometry properties of the Gauss measure, especially
the fact that the Gauss measure is doubling on each class of
admissible balls $\cb_a$ (see \cite[Proposition~2.1]{mm}), for any given $a\in(0,
\fz)$ and $\bz\in(0, 1)$, we define the local fractional integral
operator $I_a^\bz$ by that for all functions $f\in L_c^\fz(\gz)$
and $x\in\rn$,
\begin{equation}\label{1.3}
    I_a^\bz(f)(x)\equiv\dint_{B(x,\, am(x))}\dfrac{f(y)}{[V(x,
    y)]^{1-\bz}}\,d\gz(y),
\end{equation}
where and in what follows, $V(x, y)\equiv\gz(B(x, |x-y|))$
and $L_c^\fz(\gz)$ denotes the set of all functions in $L^\fz(\gz)$
with compact support. In fact, obviously, $L^\fz(\gz)=
L^\fz(\rn)$ and $L_c^\fz(\gz)=L_c^\fz(\rn)$ with equivalent norms.
To characterize the space $\BMO(\gz)$, we also introduce a variant of
the above local fractional integral operator $I_a^\bz$,
which we denote by ${\wz I}_a^\bz$.
Precisely, for all functions $f\in L_c^\fz(\gz)$ and $x\in\rn$,
\begin{equation}\label{1.3'}
    \wz{I}_a^\bz(f)(x)\equiv\dint_{B(x,\, am(x))}
    \dfrac{f(y)}{[e^{-|x|^2}|x-y|^n]^{1-\bz}}\,d\gz(y).
\end{equation}
By \eqref{2.3} below, it is not difficult to see that
there exists a positive constant $C$,
depending only on $n$ and $a$, such that for all $x\in\rn$ and $y\in B(x, am(x))$,
$$C^{-1}V(x, y)\le e^{-|x|^2}|x-y|^n\le CV(x, y).$$
Hence, when $f$ is a non-negative function,
$I_a^\bz(f)$ and $\wz{I}_a^\bz(f)$ are pointwisely equivalent.

It is proved in Theorems \ref{t3.1} and \ref{t3.2} below that
both $I_a^\bz$ and $\wz{I}_a^\bz$ are bounded from $L^p(\gz)$ to $L^{p/(1-p\bz)}(\gz)$
when $p\in(1, 1/\bz)$, or from $H^1(\gz)$ to $L^{1/(1-\bz)}(\gz)$,
or from  $L^{1/\bz}(\gz)$ to $\BMO(\gz)$
(actually, from $\{f\in L^{1/\bz}(\gz):\,f\ge0\}$ to $\blo_a(\gz)$, where $\blo_a(\gz)$
was introduced in \cite{ly} and $\blo_a(\gz)\subsetneq\BMO(\gz)$).
These results of boundedness are of independent interest;
see \cite{a75, sw60} for the corresponding  boundedness results
of the classical fractional integral operators.

If $b\in\BMO(\gz)$, then the commutator
$[b, \wz{I}_a^\bz]$, generated by $b$ and the
local fractional integral operator $\wz{I}_a^\bz$, is defined by
setting, for all functions $f\in L_c^\fz(\gz)$,
\begin{equation*}
\big[b, \wz{I}_a^\bz\big](f)\equiv b\wz{I}_a^\bz(f)-\wz{I}_a^\bz(bf).
\end{equation*}
Moreover, we define $\wz{[b, I_a^\bz]}$ by setting, for all functions $f\in L_c^\fz(\gz)$
and all $x\in\rn$,
\begin{equation*}
\wz{\big[b, I_a^\bz\big]}(f)(x)\equiv \dint_{B(x,\, am(x))}\dfrac{|b(x)-b(y)||f(y)|}{[V(x,
    y)]^{1-\bz}}\,d\gz(y).
\end{equation*}

Applying the boundedness of $I_a^\bz$ and $\wz{I}_a^\bz$
from $L^p(\gz)$ to $L^{p/(1-p\bz)}(\gz)$
with $p\in(1, 1/\bz)$,
we characterize the space $\BMO(\gz)$
by these commutators as follows.

\begin{thm}\label{t1.1}
Let $a\in(0, \fz)$, $\bz\in(0, 1)$, $1<p<q<\fz$ and $1/q=1/p-\bz$.
Then there exists a positive constant $C$, depending on $a$, $p$ and $q$,
such that the following hold.
\begin{enumerate}
\vspace{-0.25cm}
\item[(i)] If $b\in\BMO(\gz)$, then for all
$f\in L_c^\fz(\gz)$,
$$\lf\|\wz{\big[b, I_a^\bz\big]}(f)\r\|_{L^q(\gz)}\le C\|b\|_\ast\|f\|_{L^p(\gz)}.$$
Moreover, the sublinear operator $\wz{[b, I_a^\bz]}$ admits a unique
bounded extension from $L^p(\gz)$ to $L^q(\gz)$ with norm at most a constant
multiple of $\|b\|_\ast$.
\vspace{-0.25cm}
\item[(ii)] If $b\in L^1(\gz)$ and $[b, \wz{I}_a^\bz]$
is bounded from $L^p(\gz)\cap L_c^\fz(\gz)$ to $L^q(\gz)$,
then $b\in\BMO(\gz)$ and
$$\|b\|_{\BMO(\gz)}\le\|b\|_{L^1(\gz)}
+C\lf\|\big[b, \wz{I}_a^\bz\big]\r\|_{L^p(\gz)\to L^q(\gz)}.$$
\vspace{-0.8cm}
\end{enumerate}
\end{thm}

The proof of Theorem \ref{t1.1} is given in Section \ref{s4}.
Observe that there exists a positive constant $C$ such that for all $x\in\rn$,
$$|[b, \wz{I}_a^\bz](f)(x)|\le C \wz{[b, I_a^\bz]}(f)(x).$$
Thus, Theorem \ref{t1.1} actually implies that the boundedness
of either $[b, \wz{I}_a^\bz]$ or $\wz{[b, I_a^\bz]}$ characterizes $b\in\BMO(\gz)$. Moreover,
for $a\in (0,\fz)$ and $\bz\in (0,1)$,
if we define the dual operator $(I_a^\bz)^\ast$ of $I_a^\bz$ by setting,
for all functions $f\in L_c^\fz(\gz)$ and $x\in\rn$,
\begin{equation*}
    (I_a^\bz)^\ast(f)(x)\equiv\dint_{B(y,\, am(y))}\dfrac{f(y)}{[V(y,
    x)]^{1-\bz}}\,d\gz(y),
\end{equation*}
and define the dual operator $({\wz I}_a^\bz)^\ast$ of ${\wz I}_a^\bz$ by setting,
for all functions $f\in L_c^\fz(\gz)$ and $x\in\rn$,
\begin{equation*}
    (\wz I_a^\bz)^\ast(f)(x)\equiv\dint_{B(y,\, am(y))}
    \dfrac{f(y)}{[e^{-|y|^2}|x-y|^n]^{1-\bz}}\,d\gz(y),
\end{equation*}
then Theorem \ref{t1.1} and all results in Section \ref{s3}
related to $I_a^\bz$ and $\wz I_a^\bz$ are also true for $(I_a^\bz)^\ast$ and  $(\wz I_a^\bz)^\ast$;
we omit the details by similarity.

For any given $a\in(0, \fz)$ and $\bz\in(0, 1)$,
we define the local fractional maximal operator
$\cm_a^\bz$ by setting, for all locally integrable functions $f$
and all $x\in\rn$,
\begin{equation}\label{1.4}
\cm_a^\bz(f)(x)\equiv\sup_{B\in\cb_a(x)}
\dfrac1{[\gz(B)]^{1-\bz}}\dint_B |f(y)|\,d\gz(y).
\end{equation}
The boundedness results for $\cm_a^\bz$
are presented in Corollary \ref{c3.1} below.

If $b\in\BMO(\gz)$, then the commutator $[b, \cm_a^\bz]$,
generated by $b$ and the local fractional maximal
operator $\cm_a^\bz$, is defined by
setting, for all functions $f\in L_c^\fz(\gz)$,
\begin{equation}\label{1.5}
\big[b, \cm_a^\bz\big](f)\equiv b\cm_a^\bz(f)-\cm_a^\bz(bf).
\end{equation}
Correspondingly,  we define $\wz{[b, \cm_a^\bz]}$ by setting, for
all functions $f\in L_c^\fz(\gz)$ and all $x\in\rn$,
\begin{equation*}
\wz{\big[b, \cm_a^\bz\big]}(f)(x)\equiv\sup_{B\in\cb_a(x)}
\dfrac1{[\gz(B)]^{1-\bz}}\dint_B |b(x)-b(y)||f(y)|\,d\gz(y).
\end{equation*}

Applying Theorem \ref{t1.1}, we also obtain the following
characterization of the space $\BMO(\gz)$.

\begin{thm}\label{t1.2}
Let $a\in(0, \fz)$, $\bz\in(0, 1)$, $1<p<q<\fz$ and $1/q=1/p-\bz$.
Then there exists a positive constant $C$, depending on $a$, $p$ and $q$,
such that the following hold.
\begin{enumerate}
\vspace{-0.25cm}
\item[(i)] If $b\in\BMO(\gz)$, then for all
$f\in L_c^\fz(\gz)$,
$$\lf\|\wz{\big[b, \cm_a^\bz\big]}(f)\r\|_{L^q(\gz)}\le C\|b\|_\ast\|f\|_{L^p(\gz)}.$$
Moreover, the sublinear operator $\wz{[b, \cm_a^\bz]}$ admits a unique
bounded extension from $L^p(\gz)$ to $L^q(\gz)$ with norm at most a constant
multiple of $\|b\|_\ast$.
\vspace{-0.25cm}
\item[(ii)] If $b\in L^1(\gz)$ and  $\wz {[b, \cm_a^\bz]}$
is bounded from $L^p(\gz)\cap L_c^\fz(\gz)$ to $L^q(\gz)$,
then $b\in\BMO(\gz)$ and
$$\|b\|_{\BMO(\gz)}\le\|b\|_{L^1(\gz)}
+C\lf\|\wz {[b, \cm_a^\bz]}\r\|_{L^p(\gz)\to L^q(\gz)}.$$
\vspace{-0.8cm}
\end{enumerate}
\end{thm}

For any given function $f$ and all $x\in\rn$, we set
$f^+(x)\equiv\max\{f(x), 0\}$ and $f^-(x)\equiv -\min\{0, f(x)\}$.
Motivated by \cite{bmr}, we obtain the following result.

\begin{thm}\label{t1.3}
Let $a\in(0, \fz)$, $\bz\in(0, 1)$, $1<p<q<\fz$ and $1/q=1/p-\bz$.
Then there exists a positive constant $C$, depending on $a$, $p$ and $q$,
such that the following hold.
\begin{enumerate}
\vspace{-0.25cm}
\item[(i)] If $b\in\BMO(\gz)$ and $b^-\in L^\fz(\gz)$, then for all
$f\in L_c^\fz(\gz)$,
$$\|[b, \cm_a^\bz](f)\|_{L^q(\gz)}\le C\lf(\|b\|_\ast+\|b^-\|_{L^\fz(\gz)}\r)\|f\|_{L^p(\gz)}.$$
\vspace{-0.8cm}
\item[(ii)] If $b\in L^1(\gz)$ and  $[b, \cm_a^\bz]$
is bounded from $L^p(\gz)\cap L_c^\fz(\gz)$ to $L^q(\gz)$,
then $b\in\BMO(\gz)$, $b^-\in L^\fz(\gz)$ and
$$\|b\|_{\BMO(\gz)}+\|b^-\|_{L^\fz(\gz)}\le\|b\|_{L^1(\gz)}
+C\|[b, \cm_a^\bz]\|_{L^p(\gz)\to L^q(\gz)},$$
where
$$\|[b, \cm_a^\bz]\|_{L^p(\gz)\to L^q(\gz)}=
\sup\{\|[b, \cm_a^\bz](f)\|_{L^q(\gz)}:\,f\in L_c^\fz(\gz),\,
\|f\|_{L^p(\gz)}=1\}.$$

\end{enumerate}
\end{thm}

We point out that in Theorem \ref{t1.3} (i),
since the operator $[b, \cm_a^\bz]$ as in \eqref{1.5}
is not even sublinear, we can not naturally extend the domain of
$[b, \cm_a^\bz]$ to the whole $L^p(\gz)$.

Finally, we make the following conventions on notation. Let
$\nn\equiv\{1,\,2,\,\cdots\}$. Define
$p'$ to be the conjugate value of $p$, namely, $1/p+1/p'=1$ for
$p\in[1, \fz]$.
Denote by $\chi_E$ the characterize function
of any set $E\subset \rn$. We also denote by $C$ a
positive constant independent of main parameters involved, which may
vary at different occurrences. Constants with subscripts do not
change through the whole paper. We use $f\ls g$  to denote $f\le
Cg$. If $f\ls g\ls f$, we then write $f\sim g$.
For any given ``normed" spaces $\cx$ and $\cy$, an operator $T$
is said to be bounded from $\cx$ to $\cy$ means that there exists
a nonnegative constant $C$ such that for all $f\in\cx$ and $T(f)\in\cy$,
$\|T(f)\|_{\cy}\le C\|f\|_\cx$; moreover, we denote by $\|T\|_{\cx\to\cy}$
the operator norm of $T$. For any $b\in L_\loc^1(\gz)$ and any ball
$B\subset\rn$, set $b_B\equiv\frac1{\gz(B)}\int_B b(y)\,d\gz(y)$.

\section{Preliminaries}\label{s2}
\hskip\parindent
Some geometry properties concerned with the Gauss measure are used
throughout the whole paper. An important one, among others, is that
the Gauss measure is indeed doubling on all balls in $\cb_a$.
Precisely, for all $\tau$, $a\in(0, \fz)$ and $B\in\cb_a$, denote by
$B_\tau^\ast$ the union of all balls $B^\pz$ that intersect $B$ such
that $r_{B'}\le\tau r_B$. It was proved in
\cite[Proposition~2.1]{mm} that
\begin{equation}\label{2.1}
\sz_{a,\,\tau}^\ast
\equiv\sup_{B\in\cb_a}\frac{\gz(B_\tau^\ast)}{\gz(B)}
\le(2\tau+1)^ne^{4a(\tau+1)+a^2},
\end{equation}
which is deduced from the property that for all $B\in\cb_a$ and
$x\in B$,
\begin{equation}\label{2.2}
e^{-2a-a^2}\le e^{|c_B|^2-|x|^2}\le e^{2a}.
\end{equation}
Moreover, it follows from \eqref{2.2} that  for all $B\in\cb_{a}$,
\begin{equation}\label{2.3}
\gz(B)\sim e^{-|c_B|^2}|B|
\end{equation}
with constant depending only on $a$ and $n$, where $|B|$ denotes the Lebesgue measure of $B$.

Recall that it was proved in \cite[(3.4)]{mm} that for all $y\in B$
and $B\in\cb_1$, $m(y)\le2m(c_B)$. An argument similar to that also
yields that for all $B\in\cb_a$ and $y\in B$,
\begin{equation}\label{2.4}
(a+1)^{-1}m(y)\le m(c_B)\le (a+1)m(y);
\end{equation}
see also \cite{ly}. Using \eqref{2.4} and following a procedure
similar  to that in \cite[Lemma~2.4]{gmst99}, we obtain the
following geometry covering lemma. Here we omit the details.

\begin{lem}\label{l2.1}
Let $\kz=\frac1{20}$ and $\{B(x_j, \kz m(x_j)/4)\}_{j\in\nn}$ be a maximal
family of disjoint balls contained in $\rn$. Set $B_j\equiv
B(x_j, \kz m(x_j))$. Then $\rn=\cup_{j\in\nn} B_j$. Moreover, for
any given $\tau\in(0, \fz)$, there exists a positive constant $N$
depending only on $n$ and $\tau$ such that for all $x\in\rn$,
$\sum_{j\in\nn}\chi_{_{(\tau B_j)}}(x)\le N$.
\end{lem}

For any given $a\in(0, \fz)$ and locally
integrable function $f$, Mauceri and Meda \cite[p.\,286]{mm} introduced
the noncentered local Hardy-Littlewood maximal function $\cm_af$ by
setting, for all $x\in\rn$,
\begin{equation}\label{2.5}
\cm_af(x)\equiv\sup_{B\in\cb_a(x)}\dfrac1{\gz(B)}\dint_B|f(y)|\,d\gz(y),
\end{equation}
and they also introduced the local sharp maximal function $f^{\sharp}$
of $f$ by setting, for all $x\in\rn$,
\begin{equation}\label{2.6}
f^{\sharp}(x)\equiv
\sup_{B\in\cb_1(x)}\dfrac1{\gz(B)}\dint_B|f(y)-f_B|\,d\gz(y).
\end{equation}
For any fixed $a\in(0, \fz)$, it is known that
$\cm_a$ is bounded on $L^p(\gz)$ when $p\in(1, \fz]$
and from $L^1(\gz)$ to $L^{1,\,\fz}(\gz)$; see \cite[Theorem~3.1]{mm}.
When $p\in[1, \fz)$, it is also proved in \cite[Theorem~3.5]{mm} that
there exists a positive
constant $C$ such that for all $f\in L^p(\gz)$,
\begin{equation}\label{2.7}
\|f\|_{L^p(\gz)}\le C\lf\{\|f\|_{L^1(\gz)}+\|f^\sharp\|_{L^p(\gz)}\r\}.
\end{equation}

For any  $a\in(0, \fz)$, $p\in[1, \fz)$ and  locally integrable
function $f$, if we set
$$ \|f\|_\ast^{\cb_a,\, p}\equiv\sup_{B\in\cb_a}
 \lf\{\dfrac1{\gz(B)}\dint_B|f(x)-f_B|^p\,d\gz(x)\r\}^{1/p},$$
then by \cite[Proposition~2.4]{mm} and \cite[Section~4]{mm},
it  is easy to see that there exists a
positive constant $C$, depending on $a$, $n$ and $p$, but not on $f$,
such that
\begin{equation}\label{2.8}
C^{-1}\|f\|_\ast\le\|f\|_\ast^{\cb_a,\, p} \le C\|f\|_\ast.
\end{equation}

Finally, we conclude this section by recalling the $\blo$-type spaces in \cite{ly}.
Let $a\in(0, \fz)$. The space $\blo_a(\gz)$ is defined
to be the set of all locally integrable functions $f$
satisfying
\begin{equation}\label{2.9}
\|f\|_{\blo_a(\gz)}\equiv\|f\|_{L^1(\gz)}+\sup_{B\in\cb_a}
\lf[\frac1{\gz(B)}\int_B f(y)\,d\gz(y)-\mathop\essinf_{x\in B} f(x)\r]<\fz.
\end{equation}
For any fixed $a\in(0, \fz)$, it was proved in \cite{ly} that
$\blo_a(\gz)\subset \BMO(\gz)$ and
there exists a positive constant $C$ depending only on $a$ and $n$
such that for all $f\in\blo_a(\gz)$,
\begin{equation}\label{2.10}
\|f\|_{\BMO(\gz)}\le C\|f\|_{\blo_a(\gz)},
\end{equation}
and that the inclusion
$\blo_a(\gz)\subset \BMO(\gz)$ is proper.

\section{Local fractional integral and maximal operators}\label{s3}

\hskip\parindent In this section,
we establish the boundedness of $I_a^\bz$, $\wz I_a^\bz$ and $\cm_a^\bz$ on various spaces.
For the Euclidean case of Theorem \ref{t3.1}
and a variant of Theorem \ref{t3.1} on general metric measure spaces,
we refer the reader to \cite[Theorem~1]{he} and \cite[Theorem~5.3]{hk},
respectively.

\begin{thm}\label{t3.1}
Let $a\in(0, \fz)$, $\bz\in(0, 1)$, $1<p<q<\fz$ and $1/q=1/p-\bz$.
Then both $I_a^\bz$ and $\wz I_a^\bz$ are bounded from
$L^p(\gz)$ to $L^q(\gz)$ and  from $L^1(\gz)$ to $L^{1/(1-\bz),\,\fz}(\gz)$.
\end{thm}

\begin{pf}\rm
Observing that $I_a^\bz(|f|)\sim\wz I_a^\bz(|f|)$, we only need to
prove the boundedness results for the operator $I_a^\bz$.
To this end, fix $r\in(0, a)$. For any $x\in\rn$, we write
\begin{eqnarray*}
    |I_a^\bz(f)(x)|&&\le\dint_{|x-y|<rm(x)}\dfrac{|f(y)|}{[V(x, y)]^{1-\bz}}\,d\gz(y)
    +\dint_{rm(x)\le|x-y|<am(x)}\cdots\equiv {\rm Z}_1+{\rm Z}_2.
\end{eqnarray*}
By \eqref{2.2} and \eqref{2.3}, we obtain
\begin{eqnarray*}
{\rm Z}_1 &&\le\sum_{j=0}^\fz
\dint_{2^{-j-1}rm(x)\le|x-y|<2^{-j}rm(x)}
\dfrac{|f(y)|}{[\gz(B(x, 2^{-j-1}r m(x)))]^{1-\bz}}\,d\gz(y)\\
&&\le\sz_{a,\,2}^\ast\sum_{j=0}^\fz[\gz(B(x, 2^{-j-1}r m(x)))]^{\bz}
\dfrac1{{\gz(B(x, 2^{-j}r m(x)))}}
\dint_{B(x,\, 2^{-j}rm(x))}|f(y)|\,d\gz(y)\\
&&\ls e^{-|x|^2\bz}[r m(x)]^{\bz n}\cm_a(f)(x).
\end{eqnarray*}
By H\"older's inequality, we have
\begin{eqnarray*}
{\rm Z}_2 &&\le\|f\|_{L^p(\gz)}\lf\{\dint_{rm(x)\le|x-y|<am(x)}
\lf(\dfrac{1}{[\gz(B(x, |x-y|))]^{1-\bz}}\r)^{p'}\,d\gz(y)\r\}^{1/{p'}}\\
&&\le\|f\|_{L^p(\gz)}
\lf\{\sum_{j=0}^{\lfloor\log_2(\frac{a}{r})\rfloor}
\dint_{2^jrm(x)\le|x-y|<2^{j+1}rm(x)} \dfrac{1} {[\gz(B(x,
2^{j}r m(x)))]^{p'(1-\bz)}}\,d\gz(y)\r\}^{1/{p'}},
\end{eqnarray*}
where $\lfloor a \rfloor$ for any $a\in\rr$ denotes the maximal integer
no more than $a$.
Notice that $B(x, 2^{j+1}rm(x))\in\cb_{2a}$ for all $0\le j\le
\lfloor\log_2(\frac{a}{r})\rfloor$. By this, \eqref{2.2},
\eqref{2.3} and $1/q=1/p-\bz$, we further obtain
\begin{eqnarray*}
{\rm Z}_2 &&\ls\|f\|_{L^p(\gz)}
\lf\{\sum_{j=0}^{\lfloor\log_2(\frac{a}{r})\rfloor}
2^{jn(1+\bz p'-p')}[r m(x)]^{n(1+\bz p'-p')}e^{-|x|^2(1+\bz p'-p')}\r\}^{1/{p'}}\\
&&\ls \|f\|_{L^p(\gz)}[r m(x)]^{n(\bz -1/p)}e^{-|x|^2(\bz-1/p)}.
\end{eqnarray*}
Notice that $$e^{-|x|^2\bz}[rm(x)]^{ n\bz }\cm_a(f)(x)
\le\|f\|_{L^p(\gz)}[rm(x)]^{n(\bz -1/p)}e^{-|x|^2(\bz-1/p)}$$ if and
only if
\begin{equation}\label{3.1}
  r\le\dfrac1{m(x)}\lf(\dfrac{\|f\|_{L^p(\gz)}e^{|x|^2/p}}
  {\cm_a(f)(x)}\r)^{p/n}.
\end{equation}
Denote by $A_1$ the set $x\in\rn$ such that the right-hand side of
\eqref{3.1} is smaller than $a$. Set $A_2\equiv\rn\setminus A_1$. If
$x\in A_1$, we take $r$ equals to the right-hand side of
\eqref{3.1}, and then obtain
\begin{equation}\label{3.2}
|I_a^\bz(f)(x)|\ls\|f\|_{L^p(\gz)}^{\bz p} \lf[\cm_a(f)(x)\r]^{1-\bz
p}.
\end{equation}
If $x\in A_2$, we then take $r=a/2$. In this case,
$\|f\|_{L^p(\gz)}e^{|x|^2/p}\ge[am(x)]^{n/p}\cm_a(f)(x)$ for all
$x\in A_2$. This combined with the estimates of ${\rm Z}_1$ and
${\rm Z}_2$ yields that for all $x\in A_2$,
\begin{equation}\label{3.3}
|I_a^\bz(f)(x)|\ls \|f\|_{L^p(\gz)}^{\bz p}
\lf[\cm_a(f)(x)\r]^{1-\bz p}.
\end{equation}
Then by \eqref{3.2}, \eqref{3.3}, $q(1-p\bz)=p$ and the boundedness
of $\cm_a$ on $L^p(\gz)$, we obtain
\begin{eqnarray*}
\|I_a^\bz(f)\|_{L^q(\gz)}^q
&&\ls \|f\|_{L^p(\gz)}^{\bz pq}
\dint_{\rn}\lf[\cm_a(f)(x)\r]^{q(1-\bz p)}\,d\gz(x)= \|\cm_a\|_{L^p(\gz)\to
L^p(\gz)}^p\|f\|_{L^p(\gz)}^q,
\end{eqnarray*}
which implies that $I_a^\bz$ is bounded from $L^p(\gz)$ to $L^q(\gz)$.

Again using  \eqref{3.2}, \eqref{3.3} and the weak $(1, 1)$ property
of $\cm_a$, we obtain that there exists a positive constant $C$ such
that for all $\lz>0$,
\begin{eqnarray*}
\gz(\{x\in\rn:\, |I_a^\bz(f)(x)|>\lz\})
&&=\gz\lf(\lf\{x\in\rn:\, C\|f\|_{L^1(\gz)}^{\bz}
\lf[\cm_a(f)(x)\r]^{1-\bz}>\lz\r\}\r)\\
&&\ls\|\cm_a\|_{L^1(\gz)\to L^{1,\,\fz}(\gz)}
\lf(\dfrac{\|f\|_{L^1(\gz)}}{\lz}\r)^{1/(1-\bz)},
\end{eqnarray*}
which implies that $I_a^\bz$ is bounded from $L^1(\gz)$ to $L^{1/(1-\bz),\,\fz}(\gz)$.
This finishes the proof of Theorem \ref{t3.1}.
\end{pf}

Applying Theorem \ref{t3.1}, we easily deduce the following conclusion.

\begin{cor}\label{c3.1}
Let $a\in(0, \fz)$, $\bz\in(0, 1)$, $1<p<q<\fz$ and $1/q=1/p-\bz$.
Then
\begin{enumerate}
\vspace{-0.2cm}
\item[(i)] $\cm_a^\bz$  is bounded from
$L^p(\gz)$ to $L^q(\gz)$ and  from $L^1(\gz)$ to
$L^{1/(1-\bz),\,\fz}(\gz)$; \vspace{-0.2cm}
\item[(ii)] $\cm_a^\bz$  is bounded from $L^{1/\bz}(\gz)$ to $L^\fz(\gz)$.
\vspace{-0.2cm}
\end{enumerate}
\end{cor}

\begin{pf}\rm
For all $x\in\rn$ and $B\in\cb_a(x)$, by \eqref{2.4},
we obtain $B\subset B(x, 2a(a+1)m(x))$.
By this, \eqref{2.2} and \eqref{2.3}, we further
obtain
\begin{eqnarray*}
\cm_a^\bz(f)(x) &&\sim\sup_{B\in\cb_a(x)}
\dfrac1{[e^{-|c_B|^2}r_B^n]^{1-\bz}}\dint_B |f(y)|\,d\gz(y)\\
&&\ls \sup_{B\in\cb_a(x)}
\dint_{B}\dfrac{|f(y)|}{[e^{-|x|^2}|x-y|^n]^{1-\bz}}\,d\gz(y)\\
&&\ls\dint_{B(x,\, 2a(a+1)m(x))}\dfrac{|f(y)|}{[V(x,
y)]^{1-\bz}}\,d\gz(y) =I_{2a(a+1)}^\bz(f)(x),
\end{eqnarray*}
which together with Theorem \ref{t3.1} implies (i). Property (ii)
follows from H\"older's inequality. This finishes the proof of
Corollary \ref{c3.1}.
\end{pf}

Applying Theorem \ref{t3.1}, we now consider the end-point
boundedness of $I_a^\bz$.

\begin{thm}\label{t3.2}
Let $a\in(0, \fz)$, $\bz\in(0, 1)$, $1<p<q<\fz$ and $1/q=1/p-\bz$.
Then
\begin{enumerate}
\vspace{-0.25cm}
\item[(i)] both $I_a^\bz$ and $\wz I_a^\bz$ are bounded from $H^1(\gz)$ to $L^{1/(1-\bz)}(\gz)$;
\vspace{-0.25cm}
\item[(ii)] both $I_a^\bz$ and $\wz I_a^\bz$ are bounded from $\{f\in L^{1/\bz}(\gz):\, f\ge0\}$
to $\blo_a(\gz)$;
\vspace{-0.25cm}
\item[(iii)] both $I_a^\bz$ and $\wz I_a^\bz$ are
bounded from $L^{1/\bz}(\gz)$ to $\BMO(\gz)$.
\vspace{-0.25cm}
\end{enumerate}
\end{thm}

\begin{pf}\rm
We only prove the results of this theorem for the operator $I_a^\bz$, since
the proof for $\wz{I}_a^\bz$ is similar but simpler.

To prove (i), first assume that $g$ is a $(1, 2/(1+\bz))$ atom supported on $B\in\cb_1$
and show that $\|I_a^\bz(g)\|_{L^{1/(1-\bz)}(\gz)}\ls1$. To this end, we write
\begin{eqnarray}\label{3.4}
\dint_{\rn}|I_a^\bz(g)(x)|^{1/(1-\bz)}\,d\gz(x)
&&=\dint_{2B}|I_a^\bz(g)(x)|^{1/(1-\bz)}\,d\gz(x)+
\dint_{(2B)^\complement}\cdots
\equiv {\rm Y}_1+{\rm Y}_2,
\end{eqnarray}
where and in what follows, $E^\complement\equiv \rn\setminus E$
for any set $E\subset \rn$.

By H\"older's inequality and \eqref{2.1} together with the fact $I_a^\bz$ is bounded from
$L^{2/(1+\bz)}(\gz)$ to $L^{2/(1-\bz)}(\gz)$ (see Theorem \ref{t3.1}), we obtain
\begin{eqnarray*}
{\rm Y}_1
&&\le\lf[\dint_{2B}|I_a^\bz(g)(x)|^{2/(1-\bz)}\,d\gz(x)\r]^{1/2}[\gz(2B)]^{1/2}\\
&&\le [\sz_{1, 2}^\ast]^{1/2} \|I_a^\bz\|_{L^{2/(1+\bz)}(\gz)
\to L^{2/(1-\bz)}(\gz)}^{1/(1-\bz)}
\lf[\dint_{B}|g(x)|^{2/(1+\bz)}\,d\gz(x)\r]^{(1+\bz)/(2-2\bz)}[\gz(B)]^{1/2}\ls1.
\end{eqnarray*}

Now we estimate ${\rm Y}_2$. For any $x\notin 2B$ satisfying
$I_a^\bz(g)(x)\neq0$,
by \eqref{1.3} and $\supp g\subset B$,
there exists $w\in B\cap B(x, am(x))$,
which combined with \eqref{2.4} and $B\in\cb_1$ yields
\begin{equation}\label{3.5}
(2a+2)^{-1}m(x)\le m(c_B)\le(2a+2)m(x).
\end{equation}
For such an $x$, by \eqref{3.5} and the triangular inequality together with $r_B\le m(c_B)$,
we obtain $B\subset B(x, (5a+4)m(x))$
and $x\in B(c_B, a^\ast m(c_B))$ with $a^\ast\equiv 2a(a+1)+1$.
Thus,
\begin{equation}\label{3.6}
{\rm Y}_2=\dint_{B(c_B,\, a^\ast m(c_B))\cap(2B)^\complement}
\lf|I_a^\gz(g)(x)\r|^{1/(1-\bz)}\,d\gz(x).
\end{equation}
For all $x\in B(c_B, a^\ast m(c_B))\cap(2B)^\complement$ satisfying $I_a^\bz(g)(x)\neq0$,
by the facts $\int_B g(x)\,d\gz(x)=0$ and  $B\subset B(x, (5a+4)m(x))$, we write
\begin{eqnarray*}
|I_a^\bz(g)(x)|&&\le\lf|\dint_{B(x,\, (5a+4)m(x))}g(y)\lf[\dfrac{1}{[V(x, y)]^{1-\bz}}-
\dfrac{1}{[V(x, c_B)]^{1-\bz}}\r]\,d\gz(y)\r|\\
&&\hs+\lf|\dint_{B(x,\, (5a+4)m(x))\setminus B(x,\, am(x))}
\dfrac{g(y)}{[V(x, y)]^{1-\bz}}\,d\gz(y)\r|\equiv {\rm I}_1+{\rm I}_2.
\end{eqnarray*}
Since  $x\in B(c_B, a^\ast m(c_B))$, then by \eqref{2.2},
\begin{equation}\label{3.7}
e^{-|c_B|^2}\sim e^{-|x|^2}.
\end{equation}
This combined with \eqref{3.5} yields
$${\rm I}_2\le \dfrac{\|g\|_{L^1(\gz)}}{[\gz(B(x, am(x)))]^{1-\bz}}
\ls \lf(e^{-|c_B|^2}[m(c_B)]^n\r)^{\bz-1}.$$

To estimate ${\rm I}_1$, notice that for any $i\in\{1,\cdots,n\}$,
a direct calculation yields that the $i$-th partial
derivative of $V(x, y)$,
\begin{equation}\label{3.8}
\dfrac{\partial V(x, y)}{\partial y_i}=\pi^{-n/2}\dfrac{x_i-y_i}{|x-y|^{2-n}}
\dint_{S^{n-1}}e^{-||x-y|\xi'+x|^2}\,d\sz(\xi'),
\end{equation}
where and in what follows,
$S^{n-1}\equiv\{x\in\rn:\, |x|=1\}$ and $d\sz$ denotes the Lebesgue
measure on the unite sphere $S^{n-1}$.

Applying \eqref{3.8} and the mean value theorem, we have
\begin{eqnarray*}
{\rm J}&&\equiv
\lf|\dfrac{1}{[V(x, y)]^{1-\bz}}-
\dfrac{1}{[V(x, c_B)]^{1-\bz}}\r|\\
&&\le\dfrac{(1-\bz)\pi^{-n/2}|y-c_B|}{[V(x, y)]^{1-\bz}[V(x, c_B)]^{1-\bz}}
\sup_{\gfz{\tz\in(0, 1)}{z=\tz y +(1-\tz)c_B}}
[V(x, z)]^{-\bz}|x-z|^{n-1}
\dint_{S^{n-1}}e^{-||x-z|\xi'+x|^2}\,d\sz(\xi').
\end{eqnarray*}
For any $x\notin 2B$, $y\in B\subset B(x, (5a+4)m(x))$
and $z=\tz y+(1-\tz)c_B$ with $\tz\in(0, 1)$,
by \eqref{2.3}, we have
$|x-y|\sim |x-c_B|\sim |x-z|$ and
\begin{equation}\label{3.9}
V(x, z)\sim V(x, y)\sim V(x, c_B)\sim e^{-|c_B|^2}|x-c_B|^n.
\end{equation}
Moreover, by \eqref{3.5}, the facts that $|c_B|m(c_B)\le1$ and that
$w\in B\cap B(x, am(x))$, we have
\begin{eqnarray}\label{3.10}
|x||x-z|&&\le(|x-w|+|w-c_B|+|c_B|)(|x-w|+|w-c_B|+|c_B-z|)\\
&&\le (am(x)+r_B+|c_B|)(am(x)+2r_B)\ls1.\noz
\end{eqnarray}
Combining \eqref{3.7}, \eqref{3.9} and \eqref{3.10} yields
\begin{eqnarray*}
{\rm J}&&\ls\dfrac{|y-c_B| |x-c_B|^{n-1} e^{-|x|^2}}{[e^{-|c_B|^2}|x-c_B|^n]^{2-\bz}}
\ls\dfrac{|y-c_B|e^{|c_B|^2(1-\bz)}}{|x-c_B|^{n(1-\bz)+1}}.
\end{eqnarray*}
By the estimate of ${\rm J}$, H\"older's inequality, \eqref{2.2} and \eqref{2.3}, we obtain
\begin{eqnarray*}
{\rm I}_1&&\ls \dfrac{e^{|c_B|^2(1-\bz)}}{|x-c_B|^{n(1-\bz)+1}}
\dint_B|y-c_B||g(y)|\,d\gz(y)\\
&&\ls\dfrac{e^{|c_B|^2(1-\bz)}}{|x-c_B|^{n(1-\bz)+1}}\|g\|_{L^{2/(1+\bz)}(\gz)}
\lf\{\dint_B|y-c_B|^{2/(1-\bz)}\,d\gz(y)\r\}^{(1-\bz)/2}\\
&&\ls\dfrac{e^{|c_B|^2(1-\bz)}}{|x-c_B|^{n(1-\bz)+1}}
[\gz(B)]^{(\bz-1)/2} r_B^{1+n(1-\bz)/2}e^{-|c_B|^2(1-\bz)/2}\ls r_B\dfrac{e^{|c_B|^2(1-\bz)}}{|x-c_B|^{n(1-\bz)+1}}.
\end{eqnarray*}
Inserting the estimates of ${\rm I}_1$ and ${\rm I}_2$ into \eqref{3.6} yields
$${\rm Y}_2\ls \dint_{B(c_B,\, a^\ast m(c_B))\cap(2B)^\complement}
\lf[r_B^{1/(1-\bz)}\dfrac{e^{|c_B|^2}}{|x-c_B|^{n+1/(1-\bz)}}+e^{|c_B|^2}m(c_B)^{-n}
\r]\,d\gz(x)\ls1.$$
The estimates of ${\rm Y}_1$ and ${\rm Y}_2$ together with \eqref{3.4} imply
$\|I_a^\bz(g)\|_{L^{1/(1-\bz)}(\gz)}\ls1$.

If $g$ is the constant function $1$, then by \eqref{2.2} and \eqref{2.3}, we have
$$|I_a^\bz(g)(x)|\sim\dint_{B(x,\, am(x))}\dfrac1{[e^{-|x|^2} |x-y|^n]^{1-\bz}}\,d\gz(y)
\ls \dint_0^{a} s^{n-1-n(1-\bz)}\,dr\ls1,$$
which together with $\gz(\rn)=1$ further implies that
 $\|I_a^\bz(g)\|_{L^{1/(1-\bz)}(\gz)}\ls1$.
Therefore, $I_a^\bz$ maps all $(1, 2/(1+\bz))$ atoms
into uniformly bounded elements of $L^{1/(1-\bz)}(\gz)$.

For any $f\in H^1(\gz)$, we can write
$f=\sum_{j=1}^\fz\lz_ja_j$, where $\{a_j\}_{j=1}^\fz$ are $(1, 2/(1+\bz))$ atoms
and $\sum_{j=1}^\fz|\lz_j|\sim\|f\|_{H^1(\gz)}$.
By the facts $H^1(\gz)\subset L^1(\gz)$
and the boundedness of $I_a^\bz$ from $L^1(\gz)$ to $L^{1/(1-\bz),\,\fz}(\gz)$,
we obtain that $I_a^\bz(f)$ is a well-defined $L^{1/(1-\bz),\,\fz}(\gz)$ function.
Moreover, for almost all $x\in\rn$,
\begin{equation}\label{3.11}
I_a^\bz(f)(x)=\sum_{j=1}^\fz\lz_jI_a^\bz(a_j)(x);
\end{equation}
see, for example, the proof of (6.7.9) in \cite[p.\,95]{g}.
Then taking $L^{1/(1-\bz)}(\gz)$ norms on both sides of \eqref{3.11} yields that
$I_a^\bz$ is bounded from $H^1(\gz)$ to $L^{1/(1-\bz)}(\gz)$. Hence, (i) holds.

Now we show (ii). For all $f\in L^{1/\bz}(\gz)$,
by \eqref{2.2} through \eqref{2.4} together with H\"older's inequality, we obtain
\begin{eqnarray}\label{3.12}
\|I_a^\bz(f)\|_{L^1(\gz)}
&&=\lf|\dint_\rn\dint_{B(x,\, am(x))}
\dfrac{f(y)}{[V(x, y)]^{1-\bz}}\,d\gz(y)\,d\gz(x)\r|\\
&&\ls \dint_\rn\lf[\dint_{|x-y|<a(a+1)m(y)}
\dfrac{e^{-|x|^2}}{[e^{-|x|^2}|x-y|^n]^{1-\bz}}\,dx\r]|f(y)|\,d\gz(y)\noz\\
&&\ls\|f\|_{L^1(\gz)}\ls\|f\|_{L^{1/\bz}(\gz)}.\noz
\end{eqnarray}
Thus, to finish the proof of  (ii), by \eqref{3.12} and \eqref{2.9}, it suffices to show that
there exists a positive constant $C$ such that for all $f\ge0$ and $B\in\cb_a$,
\begin{equation}\label{3.13}
\frac1{\gz(B)}\int_B I_a^\bz(f)(y)\,d\gz(y)\le C\|f\|_{L^{1/\bz}(\gz)}+
\mathop\essinf_{x\in B} I_a^\bz(f)(x).
\end{equation}
To see \eqref{3.13}, we decompose $f$ into $f=f\chi_{3B}+f\chi_{(3B)^\complement}$.
Choose $p\in(1, 1/\bz)$ and $1/q=1/p-\bz$.
Using H\"older's inequality and \eqref{2.1} together with Theorem \ref{t3.1}, we obtain
\begin{eqnarray*}
\frac1{\gz(B)}\int_B I_a^\bz(f\chi_{3B})(y)\,d\gz(y)
&& \le \lf\{\frac1{\gz(B)}\int_B |I_a^\bz(f\chi_{3B})(y)|^q\,d\gz(y)\r\}^{1/q}\\
&& \le \|I_a^\bz\|_{L^p(\gz)\to L^q(\gz)} [\gz(B)]^{-1/q}\|f\chi_{3B}\|_{L^{p}(\gz)}\noz\\
&&\le[\sz_{a,\,3}^\ast]^{1/q}\|I_a^\bz\|_{L^p(\gz)\to L^q(\gz)}\|f\|_{L^{1/\bz}(\gz)}.\noz
\end{eqnarray*}
By this and the linearity of $I_a^\bz$, we know that to obtain \eqref{3.13},
it is enough to show that for all $y\in B$ and $x\in B$,
\begin{equation}\label{3.14}
I_a^\bz(f\chi_{(3B)^\complement})(y)
\le C\|f\|_{L^{1/\bz}(\gz)}+
I_a^\bz(f)(x).
\end{equation}
If $I_a^\bz(f)(x)=\fz$, then \eqref{3.14} holds trivially.
Assume now that $I_a^\bz(f)(x)<\fz$.
Notice that $f\ge0$, and hence $I_a^\bz(f)(x)$ is finite.
Thus,
\begin{eqnarray*}
&&I_a^\bz(f\chi_{(3B)^\complement})(y)-I_a^\bz(f)(x)\\
&&\hs\le\dint_{B(y,\, am(y))}\dfrac{f(z)\chi_{(3B)^\complement}(z)}{[V(y,z)]^{1-\bz}}\,d\gz(z)
-\dint_{B(x,\, am(x))}\dfrac{f(z)\chi_{(3B)^\complement}(z)}{[V(x,z)]^{1-\bz}}\,d\gz(z)
\equiv {\rm Y}.
\end{eqnarray*}
If $B\in\cb_a$, $x\in B$ and  $y\in B$, then by \eqref{2.4}, we obtain that
for all $z\in B(y, am(y))$,
$$|z-x|\le|z-y|+|y-x|< am(y)+2r_B\le a(a+1)m(c_B)+2am(c_B)\le\wz a m(x),$$
where $\wz a\equiv a(a+1)(a+3)$. Thus, $B(y, am(y))\subset B(x, \wz a m(x))$.
It follows that
\begin{eqnarray*}
{\rm Y}
&&\le \dint_{B(x,\, \wz a m(x))}\dfrac{f(z)\chi_{(3B)^\complement}(z)}{[V(y,z)]^{1-\bz}}\,d\gz(z)
-\dint_{B(x,\, am(x))}\dfrac{f(z)\chi_{(3B)^\complement}(z)}{[V(x,z)]^{1-\bz}}\,d\gz(z)\\
&&=\dint_{B(x,\, \wz a m(x))}f(z)\chi_{(3B)^\complement}(z)
\lf[\dfrac{1}{[V(y,z)]^{1-\bz}}-\dfrac1{[V(x,z)]^{1-\bz}}\r]\,d\gz(z)\\
&&\hs+\dint_{B(x,\, \wz a m(x))\setminus B(x,\, am(x))}
\dfrac{f(z)\chi_{(3B)^\complement}(z)}{[V(x,z)]^{1-\bz}}\,d\gz(z)
\equiv {\rm J}_1 +{\rm J}_2.
\end{eqnarray*}
Using H\"older's inequality and \eqref{2.1} yields
$${\rm J}_2\le\dfrac{1}{[\gz(B(x, am(x)))]^{1-\bz}}
\dint_{B(x,\, \wz a m(x))}f(z)\,d\gz(z)
\ls \|f\|_{L^{1/\bz}(\gz)}.$$
Now we estimate ${\rm J}_1$. Notice that for all $i\in\{1,\cdots, n\}$,
a simple calculation yields
\begin{eqnarray*}
\dfrac{\partial V(y, z)}{\partial y_i}
&&=\pi^{-n/2}\dfrac{y_i-z_i}{|y-z|^{2-n}}
\dint_{S^{n-1}}e^{-||y-z|\xi'+z|^2}\,d\sz(\xi')\\
&&\hs-2\pi^{-n/2}\dint_{|\xi-y|<|y-z|} (\xi_i-y_i) e^{-|\xi|^2}\,d\xi.
\end{eqnarray*}
This together with an argument similar to the estimates of ${\rm J}$ and
${\rm I}_1$ implies that
\begin{eqnarray*}
{\rm J}_1&&\le \dint_{B(x,\, \wz a m(x))}
\dfrac{(1-\bz)f(z)\chi_{(3B)^\complement}(z)|y-x|}{[V(y,z)]^{1-\bz}[V(x,z)]^{1-\bz}}
\sup_{\gfz{\tz\in(0, 1)}{w=\tz x+(1-\tz) y}}\lf|[V(w, z)]^{-\bz}
\dfrac{\partial V(w, z)}{\partial w_i}\r|\,d\gz(z)\\
&&\ls r_B\dint_{B(x,\, \wz a m(x))}
\dfrac{f(z)\chi_{(3B)^\complement}(z)}{[e^{-|x|^2}|x-z|^n]^{2-\bz}}
\sup_{\gfz{\tz\in(0, 1)}{w=\tz x+(1-\tz) y}}\lf\{|w-z|^{n-1}
\dint_{S^{n-1}}e^{-||w-z|\xi'+z|^2}\,d\sz(\xi')\r.\\
&&\hs\lf.
+\dint_{|\xi-w|<|w-z|} |\xi-w| e^{-|\xi|^2}\,d\xi\r\}\,d\gz(z)\\
&&\ls r_B\dint_{B(x,\, \wz a m(x))}
\dfrac{f(z)\chi_{(3B)^\complement}(z)}{[e^{-|x|^2}|x-z|^n]^{2-\bz}}
e^{-|z|^2}|x-z|^{n-1}\,d\gz(z)\ls\|f\|_{L^{1/\bz}(\gz)},
\end{eqnarray*}
where in the penultimate inequality, we used the fact that $|w-z|\sim|x-z|\le \wz a$
and \eqref{2.2} together with $|\xi-z|\ls m(z)$.
Combining the estimates of ${\rm J}_1$ and ${\rm J}_2$ yields the desired
estimate of ${\rm Y}$, and hence \eqref{3.14}.
Thus, (ii) holds.

Finally we prove (iii). For any $f\in L^{1/\bz}(\gz)$, we decompose
$f=f^+-f^-$, where $f^+\equiv \max \{f, 0\} $ and $f^{-}\equiv \min\{f, 0\}$.
From \eqref{2.10} and (ii) of this theorem, it follows that
\begin{eqnarray*}
\|I_a^\bz(f)\|_{\BMO(\gz)}
&&\le\|I_a^\bz(f^+)\|_{\BMO(\gz)}+\|I_a^\bz(f^-)\|_{\BMO(\gz)}\\
&&\ls \|I_a^\bz(f^+)\|_{\blo_a(\gz)}+\|I_a^\bz(f^-)\|_{\blo_a(\gz)}\ls \|f\|_{L^{1/\bz}(\gz)}.
\end{eqnarray*}
Thus (iii) holds, which completes the proof of Theorem \ref{t3.2}.
\end{pf}

\section{Proofs of Theorems \ref{t1.1} through \ref{t1.3}}\label{s4}

\hskip\parindent
In this section, by applying Theorem \ref{t3.1} and
the geometry properties listed in Section \ref{s2}
together with some ideas used in the Euclidean case
(see, for example, \cite{j78, ch82, ldy, g}),
we prove Theorems \ref{t1.1} and \ref{t1.2}.

We begin with two technical lemmas.
For any $a\in(0, \fz)$, $\bz\in(0, 1)$ and $b\in\BMO(\gz)$,
define the following auxiliary operator $T_a^\bz(b;\cdot)$ by setting, for all
locally integrable functions $f$ and $x\in\rn$,
\begin{equation}\label{4.1}
T_a^\bz(b; f)(x)\equiv \dint_{B(x,\, am(x))}
 \dfrac{|b(x)-b(y)||f(y)|}{[e^{-|y|^2}|x-y|^n]^{(1-\bz)}}\,d\gz(y).
 \end{equation}
Correspondingly, we introduce another auxiliary operator
$\wz{T_a^\bz(b;\cdot)}$,
which is in fact a smooth version of $T_a^\bz(b;\cdot)$.
Precisely, let $\phi$ be a radial function in $C_c^\fz(\rr)$
such that $0\le\phi\le1$, $\phi(t)\equiv1$ when $|t|<1$,
$\phi(t)\equiv0$ when $|t|\ge2$;
moreover, there exists a positive constant $C$ such that
$|\phi^\pz(t)|\le C/|t|$ for all $t\in\rr$.
For any $x$, $y\in\rn$, we set
$$\phi_{x,\,y}\equiv\phi\lf(\dfrac{|x-y|}{a(a+1)m(y)}\r).$$
Define
\begin{equation}\label{4.x1}
\wz{T_a^\bz(b;f)}(x)\equiv\dint_\rn
 \dfrac{\phi_{x,y}|b(x)-b(y)||f(y)|}{[e^{-|y|^2}|x-y|^n]^{(1-\bz)}}\,d\gz(y).
 \end{equation}
By \eqref{2.4} and the support condition of $\phi$, it is not difficult to see that
for all locally integrable functions $f$ and $x\in \rn$,
\begin{equation}\label{4.x2}
T_a^\bz(b;f)(x)\le\wz{T_a^\bz(b;f)}(x)\le T_{\wz C_a}^\bz(b;f)(x),
\end{equation}
where $\wz C_a\equiv 2a(a+1)(2a^2+2a+1)$.
It follows from  \eqref{2.2} that for all $x\in\rn$,
\begin{equation}\label{4.2}
\lf|\big[b, \wz I_a^\bz\big](f)(x)\r|\le
\wz{\big[b, I_a^\bz\big]}(f)(x)\sim T_a^\bz(b; f)(x).\end{equation}

The following lemma is used in the proof
of Theorem \ref{t1.1} due to the extra term $\|\cdot\|_{L^1(\gz)}$
appearing in \eqref{2.7}, comparing to the classical case of Euclidean spaces
(see, for example, \cite[p.\,148]{st93}).

\begin{lem}\label{l4.1}
Let $a\in(0, \fz)$, $\bz\in(0, 1)$, $p\in(1,1/\bz)$ and  $b\in \BMO(\gz)$.
Then both $T_a^\bz(b;\cdot)$ and $\wz{[b, I_a^\bz]}$ are
bounded from $L^p(\gz)$ to $L^1(\gz)$ with norm at most
a constant multiple of $\|b\|_\ast$.
\end{lem}

\begin{pf}\rm
By \eqref{4.2}, we only need to consider $T_a^\bz(b;\cdot)$.
To this end, let $\{B_j\}_{j\in\nn}$ be the sequence of balls as in Lemma
\ref{l2.1}. If $x\in B_j$, then by \eqref{2.4}, we obtain that for
all $y\in B(x, am(x))$,
$|y-x_j|<(2a+1)m(x_j).$
Thus, we have $B(x, am(x))\subset \frac1\kz(2a+1)B_j$,
where $\kz=\frac1{20}$ is as in Lemma \ref{l2.1}.
Set $\wz a\equiv \frac1\kz(2a+1)$.
Then for any $z\in\wz aB_j$, again using \eqref{2.4}, we
obtain $|z-x|<4(a+1)m(x),$ and
hence $\wz aB_j\subset B(x, 4(a+1)m(x))$.
From these facts together with $\rn\subset\cup_{j\in\nn} B_j$, it follows that
\begin{eqnarray*}
\|T_a^\bz(b; f)\|_{L^1(\gz)} &&\ls\sum_{j\in\nn} \dint_{B_j}\dint_{B(x,\, am(x))}
\dfrac{|b(x)-b(y)||f(y)|}{[V(x, y)]^{1-\bz}}\,d\gz(y)\,d\gz(x)\\
&&\ls\sum_{j\in\nn} \dint_{B_j}\dint_{\wz aB_j}
\dfrac{|b(x)-b(y)||f(y)|}{[V(x, y)]^{1-\bz}}\,d\gz(y)\,d\gz(x)\\
&&\ls\sum_{j\in\nn}\dint_{B_j}\lf|b(x)-b_{\wz aB_j}\r|
 I_{4(a+1)}^\bz\lf(|f|\r)(x)\,d\gz(x)\\
&&\quad+ \sum_{j\in\nn}\dint_{B_j}
I_{4(a+1)}^\bz\lf(\lf|b-b_{\wz aB_j}\r||f|\chi_{\wz aB_j}\r)(x)\,d\gz(x)
\equiv {\rm I}+{\rm J}.
\end{eqnarray*}

Choose $q\in(0, \fz)$ satisfying that $1/q=1/p-\bz$. Applying
H\"older's inequality for integrals and series, Lemma \ref{l2.1}, Theorem \ref{t3.1},
\eqref{2.1}, \eqref{2.8} and the fact $\gz(\rn)=1$, we obtain
\begin{eqnarray*}
{\rm I} &&\le\sum_{j\in\nn}
\lf\{\dint_{B_j}\lf|b(x)-b_{\wz aB_j}\r|^{q'}\,d\gz(x)\r\}^{1/q'}
\lf\{\dint_{B_j}|I_{4(a+1)}^\bz\lf(|f|\r)(x)|^q\,d\gz(x)\r\}^{1/q}\\
&&\ls\lf\{\sum_{j\in\nn}\dint_{B_j}
\lf|b(x)-b_{\wz aB_j}\r|^{q'}\,d\gz(x)\r\}^{1/q'}
\lf\|I_{4(a+1)}^\bz\lf(|f|\r)\r\|_{L^q(\gz)}\ls \|b\|_\ast\,\|f\|_{L^p(\gz)}.
\end{eqnarray*}

Now we estimate ${\rm J}$. Choose $s\in(1, p)$ and $r\in(1, q)$
such that $1/r=1/s-\bz$. Then by H\"older's inequality and Theorem \ref{t3.1},
\begin{eqnarray*}
{\rm J}&&\le \sum_{j\in\nn}[\gz(B_j)]^{1/r'}
\lf\{\dint_{B_j}\lf[I_{4(a+1)}^\bz
\lf(|b-b_{\wz aB_j}||f|\chi_{\wz aB_j}\r)(x)\r]^r\,d\gz(x)\r\}^{1/r}\\
&&\ls \sum_{j\in\nn}[\gz(B_j)]^{1/r'}
\lf\{\dint_{\wz aB_j}\lf[|b(x)-b_{\wz aB_j}||f(x)|\r]^s
\,d\gz(x)\r\}^{1/s}.
\end{eqnarray*}
Using H\"older's inequality again yields that the last quantity above is bounded by
$$\sum_{j\in\nn}[\gz(B_j)]^{1/r'}
\lf\{\dint_{\wz a B_j}\lf|b(x)-b_{\wz a B_j}\r|^{ps/(p-s)}
\,d\gz(x)\r\}^{1/s-1/p}\lf\{\dint_{\wz a B_j}|f(x)|^p
\,d\gz(x)\r\}^{1/p}.$$
Then by \eqref{2.1}, H\"older's inequality for series, \eqref{2.8}, Lemma \ref{l2.1}
and $\gz(\rn)=1$, we finally obtain
\begin{eqnarray*}
{\rm J}&&\ls\|b\|_\ast\sum_{j\in\nn}[\gz(B_j)]^{1/r'+1/s-1/p}\lf\{\dint_{\wz a B_j}|f(x)|^p
\,d\gz(x)\r\}^{1/p}\\
&&\ls\|b\|_\ast\,\|f\|_{L^p(\gz)}
\lf\{\sum_{j\in\nn}[\gz(B_j)]^{p'(1/r'+1/s-1/p)}\r\}^{1/p'}
\ls\|b\|_\ast\,\|f\|_{L^p(\gz)},
\end{eqnarray*}
where in the last step we use $p'(1/r'+1/s-1/p)=1+p'\bz$
and $\sum_{j\in\nn}[\gz(B_j)]^{1+p'\bz}\ls1$. Combining the estimates
of ${\rm I}$ and ${\rm II}$ yields the desired result for
$T_a^\bz(b;\cdot)$, which completes the proof of Lemma
\ref{l4.1}.
\end{pf}

For the Euclidean version of the following Lemma \ref{l4.2},
see \cite[Lemma~11]{j78} or \cite[p.\,418]{tor}.

\begin{lem}\label{l4.2}
Let $a\in(0, \fz)$, $\bz\in(0, 1)$, $1<p<q<\fz$ and $1/q=1/p-\bz$.
Let $b\in\BMO(\gz)$ and $\wz{T_a^\bz(b; f)}$ be as in \eqref{4.x1}.
Then for any fixed $r,\,s\in(1, p)$, there exists a positive
constant $C$ such that for all $f\in L_c^\fz(\gz)$ and $x\in\rn$,
\begin{equation}\label{4.3}
\lf(\wz{T_a^\bz(b; f)}\r)^\sharp(x)\le C\|b\|_\ast\lf\{
  \lf[\cm_{1}\lf(\big[I_{C_a}^\bz(|f|)\big]^r\r)(x)\r]^{1/r}
  +\lf[\cm_{C_a}^{\bz s}(|f|^s)(x)\r]^{1/s}\r\},
\end{equation}
where $C_a$ is a sufficiently large positive constant depending
only on $a$, $\lf(\wz{T_a^\bz(b; f)}\r)^\sharp$
is as in \eqref{2.6} with $f$ replaced by $\wz{T_a^\bz(b; f)}$,
$\cm_1$ is as in \eqref{2.5} with $a=1$
and $\cm_{C_a}^{\bz s}$ is as in \eqref{1.4}
with $\bz$ and $a$ there replaced,
respectively, by $\bz s$ and $C_a$.
\end{lem}

\begin{pf}\rm Fix $x\in\rn$.
To show \eqref{4.3}, for any fixed ball $B'\in\cb_1(x)$ and any
$y\in B'$, we write
\begin{eqnarray*}
 \wz{T_a^\bz(b; f)}(y)
&&\le |b(y)-b_{B'}|\dint_\rn
 \dfrac{\phi_{y,z}|f(z)|}{[e^{-|z|^2}|y-z|^n]^{1-\bz}}\,d\gz(z)\\
&&\quad + \dint_\rn
 \dfrac{\phi_{y,z}|b_{B'}-b(z)||f(z)|\chi_{3B'}(z)}{[e^{-|z|^2}|y-z|^n]^{1-\bz}}\,d\gz(z)\\
&&\quad +\dint_\rn
 \dfrac{\phi_{y,z}|b_{B'}-b(z)||f(z)|\chi_{(3B')^\complement}(z)}
 {[e^{-|z|^2}|y-z|^n]^{1-\bz}}\,d\gz(z)\\
&&\equiv \cj_1^{B'}(y)+\cj_2^{B'}(y)+\cj_3^{B'}(y).
\end{eqnarray*}

For all $z\in\rn$ satisfying $\phi_{y,z}\neq0$, i.\,e,
$|z-y|<2a(a+1)m(z)$, by \eqref{2.2} and \eqref{2.3}, we have,
\begin{equation}\label{4.4}
e^{-|z|^2}|z-y|^n\sim e^{-|y|^2}|z-y|^n\sim V(z, y).
\end{equation}
This and the support condition of $\phi_{y,z}$ together
with H\"older's inequality and \eqref{2.8} imply that
\begin{equation}\label{4.5}
\dfrac1{\gz(B')}\dint_{B'}\cj_1^{B'}(y)\,d\gz(y)\ls\|b\|_\ast
\lf[\cm_{1}\lf(\big[I_{C_a}^\bz(|f|)\big]^r\r)(x)\r]^{1/r},
\end{equation}
for some positive constant $C_a$ which depends only on $a$.

Choose $\dz>1$ and $\kz>1$ such that $\dz\kz=s$. Since
$1<\kz<p<1/\bz$, there exists $u>\kz$ such that $1/u=1/\kz-\bz$. By
\eqref{4.4}, Theorem \ref{t3.1}, H\"older's inequality and $\dz\kz=s$,
we obtain that for sufficiently large positive number $C_a$,
\begin{eqnarray*}
\dfrac1{\gz(B')}\dint_{B'}\cj_2^{B'}(y)\,d\gz(y)
&&\ls\lf\{\dfrac1{\gz(B')}\dint_{B'}
|I_{C_a}^{\bz}(|b_{B'}-b||f|\chi_{3B'})(y)|^u\,d\gz(y)\r\}^{1/u}\noz\\
&&\le\dfrac1{[\gz(B')]^{1/u}}
\lf\{\dint_{3B'}|b(y)-b_{B'}|^{\kz\dz'}\,d\gz(y)\r\}^{1/{(\kz\dz')}}
\lf\{\dint_{3B'}|f(y)|^{s}\r\}^{1/s}.\noz
\end{eqnarray*}
Notice that the triangular inequality of $\|\cdot\|_{L^{\kz\dz'}(\gz)}$
and \eqref{2.1} imply
\begin{eqnarray*}
\lf\{\dfrac1{\gz(3B')}
\dint_{3B'}|b(y)-b_{B'}|^{\kz\dz'}\,d\gz(y)\r\}^{1/{(\kz\dz')}}
&&\le\|b\|_\ast^{\cb_3,\,\kz\dz'}+|b_{3B'}-b_{B'}|\\
&&\le\|b\|_\ast^{\cb_3,\,\kz\dz'}+\sz_{1,\,
3}^\ast\|b\|_\ast^{3,\,1}.
\end{eqnarray*}
This combined with \eqref{2.8} and \eqref{2.1}  yields
\begin{equation}\label{4.6}
\dfrac1{\gz(B')}\dint_{B'}\cj_2^{B'}(y)\,d\gz(y)\ls\|b\|_\ast
\lf[\cm_3^{\bz s}(|f|^s)(x)\r]^{1/s}.\end{equation}

Let $y\in B'$ and $z\notin 3B'$ satisfying that
$\frac{\phi_{y,z}}{|y-z|^{n(1-\bz)}}
-\frac{\phi_{c_{B'},\,z}}{|c_{B'}-z|^{n(1-\bz)}}\neq0$.
Then we use \eqref{2.4} to obtain that
$|z-y|\le C_am(y)$ for some large enough positive constant $C_a$ that depends only on $a$.
From this, the mean value theorem, the definition of $\phi$
and the fact that $|z-y|\sim|z-c_{B'}|$ for all $y\in B'$ and $z\notin 3B'$,
we easily deduce that
$$\lf|\frac{\phi_{y,z}}{|y-z|^{n(1-\bz)}}
-\frac{\phi_{c_{B'},\,z}}{|c_{B'}-z|^{n(1-\bz)}}\r|
\ls\dfrac{|y-c_{B'}|}{|y-z|^{n(1-\bz)+1}}\chi_{B(y, C_am(y))}(z),$$
which combined with H\"older's inequality further implies that for all
$y\in B'$,
\begin{eqnarray*}
|\cj_3^{B'}(y)-\cj_3^{B'}(c_{B'})|
 &&\ls\dint_{(3B')^\complement\cap B(y,\, C_a m(y))}
 \dfrac{|b_{B'}-b(z)||f(z)||y-c_{B'}|}
 {e^{-|z|^2(1-\bz)}|z-c_{B'}|^{1+n-n\bz}}\,d\gz(z)\\
 &&\le\lf\{r_{B'}\dint_{(3B')^\complement\cap B(y,\, C_a m(y))}
 \dfrac{|b_{B'}-b(z)|^{s'}}
 {e^{-|z|^2}|z-c_{B'}|^{1+n}}\,d\gz(z)\r\}^{1/s'}\\
 &&\hs\times\lf\{r_{B'}\dint_{(3B')^\complement\cap B(y,\, C_a  m(y))}
 \dfrac{|f(z)|^s}
 {e^{-|z|^2(1-\bz s)}|z-c_{B'}|^{1+n-n\bz s}}\,d\gz(z)\r\}^{1/s}\\
&&\equiv {\rm I}\times {\rm J}.
\end{eqnarray*}
For all $z\in B(y, C_a m(y))$, by \eqref{2.2} and $y\in B'$, we have
$e^{-|z|^2}\sim e^{-|y|^2}\sim e^{-|c_{B'}|^2}$. Thus,
\begin{eqnarray*}
{\rm I}\ls \sum_{j=1}^\fz \lf\{
\dfrac{2^{-j}}{\gz(2^jB')}\dint_{\gfz{|z-y|<C_a  m(y)}{2^j
r_{B'}<|z-c_{B'}|\le 2^{j+1} r_{B'}}}
 |b_{B'}-b(z)|^{s'}  \,d\gz(z)\r\}^{1/s'}.
\end{eqnarray*}
Notice that for $j\in\nn$ satisfying $2^jr_{B'}<|z-c_{B'}|$ and $|z-y|<C_a m(y)$, we have
$$2^j r_{B'}<|z-c_{B'}|\le|z-y|+|y-c_{B'}|<C_a m(y)+r_{B'}\le[C_a(a+1)+a] m(c_{B'}),$$
and hence, $2^{j+1}B'\in\cb_{C_a(a+1)+a}$.
For simplicity, we set $ a^\ast\equiv C_a(a+1)+a$.
This, together with Minkowski's inequality,
\eqref{2.1} and \eqref{2.8}, implies that
\begin{eqnarray*}
{\rm I} &&\ls \sum_{\{j\in\nn:\,2^j r_{B'}<a^\ast m(c_{B'})\}}\hs
2^{-j/s'}
\bigg(|b_{B'}-b_{2B'}|+\cdots +|b_{2^jB'}-b_{2^{j+1} B'}|\\
&&\quad+\lf\{ \dfrac{1}{\gz(2^jB')}\dint_{2^{j+1}B'}
 |b_{2^{j+1}B'}-b(z)|^{s'}  \,d\gz(z)\r\}^{1/s'}\bigg)\\
&&\ls\sum_{j=1}^\fz 2^{-j/s'}\bigg(j\|b\|_\ast^{\cb_{a^\ast},\,1}
+\|b\|_\ast^{\cb_{2a^\ast},\,s'}\bigg)\ls\|b\|_\ast .
\end{eqnarray*}
Similarly to the estimate of ${\rm I}$, we have
\begin{eqnarray*}
{\rm J}&&\ls \sum_{j=1}^\fz \lf\{
\dfrac{2^{-j}}{[\gz(2^jB')]^{1-\bz s}}\dint_{\gfz{|z-y|<a^\ast m(y)}{2^j
r_{B'}<|z-c_{B'}|\le 2^{j+1} r_{B'}}}
 |f(z)|^{s}  \,d\gz(z)\r\}^{1/s}\\
 &&\ls \sum_{\{j\in\nn:\,2^j r_{B'}<a^\ast m(c_{B'})\}}\lf\{
\dfrac{2^{-j}}{[\gz(2^jB')]^{1-\bz s}}\dint_{2^{j+1}B'}
 |f(z)|^{s}  \,d\gz(z)\r\}^{1/s}\ls [\cm_{2a^\ast}^{\bz s}(|f|^s)(x)]^{1/s}.
\end{eqnarray*}
Combining the estimates of ${\rm I}$ and {\rm J} yields
\begin{equation}\label{4.7}
   \dfrac1{\gz(B')}\dint_{B'}
   |\cj_3^{B'}(y)-\cj_3^{B'}(c_{B'})|\,d\gz(y)
   \ls\|b\|_\ast \lf[\cm_{2a^\ast}^{\bz s}(|f|^s)(x)\r]^{1/s}.
\end{equation}
Applying \eqref{4.5}, \eqref{4.6} and \eqref{4.7},
we obtain \eqref{4.3}, which completes the proof of Lemma \ref{l4.2}.
\end{pf}

\newtheorem{prf}{\bf Proof of Theorem \ref{t1.1}}
\renewcommand\theprf{}

\begin{prf}\rm
To show  (i), we let $b\in\BMO(\gz)$. For any given $N\in\nn$, set
$b_N\equiv -N\chi_{\{b<-N\}} +b\chi_{\{|b|\le N\}}+N\chi_{\{b>N\}}.$
Then $\|b_N\|_\ast\le 4\|b\|_\ast$; see \cite[pp.\,631-632]{cw2}.
Let $\wz C_a$ and $C_a$ be respectively as in \eqref{4.x2} and Lemma \ref{l4.2}.
For all $f\in L_c^\fz(\gz)$, combining \eqref{4.1}, \eqref{4.x2},
\eqref{2.2} and Theorem \ref{t3.1} yields
$$\lf\|\wz{T_a^\bz(b_N; f)}\r\|_{L^q(\gz)}
\le\lf\|T_{\wz C_a}^\bz(b_N; f)\r\|_{L^q(\gz)}
\ls N\lf\|I^\bz_{\wz C_a}(|f|)\r\|_{L^q(\gz)} \ls N\|f\|_{L^p(\gz)}<\fz.$$
This allows us to use \eqref{2.7}. Thus, by \eqref{2.7},
Lemma \ref{l4.2} together with the fact $1<r,\,s<p$,
the $L^{q/r}(\gz)$-boundedness of $\cm_1$,
Theorem \ref{t3.1}, Corollary \ref{c3.1} and Lemma \ref{l4.1}, we obtain
\begin{eqnarray*}
&&\lf\|T_a^\bz(b_N; f)\r\|_{L^q(\gz)}\\
&&\hs\le\lf\|\wz{T_a^\bz(b_N; f)}\r\|_{L^q(\gz)}\\
&&\hs\ls\lf\|\big(\wz{T_a^\bz(b_N; f)}\big)^\sharp\r\|_{L^q(\gz)}+
\lf\|\wz{T_a^\bz(b_N; f)}\r\|_{L^1(\gz)}\\
&&\hs\ls \|b_N\|_\ast\lf\{
\lf\|\lf[\cm_{1}\lf(\big[I_{C_a}^\bz(|f|)\big]^r\r)\r]^{1/r}\r\|_{L^q(\gz)}
+\lf\|[\cm_{C_a}^{\bz s}(|f|^s)]^{1/s}\r\|_{L^q(\gz)} +\|f\|_{L^p(\gz)}\r\}\\
&&\hs\ls\|b\|_\ast\,\|f\|_{L^p(\gz)}.
\end{eqnarray*}
The dominated convergence theorem gives that $b_N\to b$ in
$L^p(\gz)$ of every compact set and, in particular, in $L^p(d\gz,
\supp(f))$. Hence, $b_N|f|\to b|f|$ in $L^p(\gz)$.
From Theorem \ref{t3.1}, it follows easily that both
\begin{equation}\label{4.8}
\dint_{\gfz{B(x,\, am(x))}{b_N(x)>b_N(y)}}
 \dfrac{b_N(y)|f(y)|}{[e^{-|y|^2}|x-y|^n]^{1-\bz}}\,d\gz(y)
 \to \dint_{\gfz{B(x,\, am(x))}{b(x)>b(y)}}
 \dfrac{b(y)|f(y)|}{[e^{-|y|^2}|x-y|^n]^{1-\bz}}\,d\gz(y)
\end{equation}
and
\begin{equation}\label{4.9}
\dint_{\gfz{B(x,\, am(x))}{b_N(x)\le b_N(y)}}
 \dfrac{b_N(y)|f(y)|}{[e^{-|y|^2}|x-y|^n]^{1-\bz}}\,d\gz(y)
 \to \dint_{\gfz{B(x,\, am(x))}{b(x)\le b(y)}}
 \dfrac{b(y)|f(y)|}{[e^{-|y|^2}|x-y|^n]^{1-\bz}}\,d\gz(y)
\end{equation}
 as $N\to\fz$ in $L^q(\gz)$. The same is true for \eqref{4.8} and \eqref{4.9}
 without $b_N$ and $b$ appearing in the integrand functions.
These observations together with the Riesz lemma and \eqref{4.1} imply that there exists
 a subsequence $\{N_k\}_{k=1}^\fz\subset\nn$ such that
$T_a^\bz(b_{N_k}; f)\to T_a^\bz(b; f)$ as $k\to\fz$ almost
everywhere. Letting $k\to\fz$ and using Fatou's lemma together with
\eqref{4.2}, we obtain that for all $f\in L_c^\fz(\gz)$,
$$\lf\|\wz{\big[b, I_a^\bz\big]}(f)\r\|_{L^q(\gz)}
\ls\lf\|T_a^\bz(b; f)\r\|_{L^q(\gz)}
\le\lim_{k\to\fz}\lf\|T_a^\bz(b_{N_k}; f)\r\|_{L^q(\gz)}
\ls\|b\|_\ast\,\|f\|_{L^p(\gz)}.$$ Applying the density of
$L_c^\fz(\gz)$ in the classical Lebesgue space $L^p(\rn)$ and the
fact that $\gz(E)\ls |E|$ for all set $E\subset\rn$, one can deduce
that $L_c^\fz(\gz)$ is dense in $L^p(\gz)$. Moreover, for any
$f,\,g\in L_c^\fz(\gz)$, it is easy to deduce that
$$\bigg|\wz{[b, I_a^\bz\big]}(f+g)\bigg|
\le \bigg|\wz{[b, I_a^\bz\big]}(f)\bigg|
+\bigg|\wz{[b, I_a^\bz\big]}(g)\bigg|$$
and
$$\lf|\wz{[b, I_a^\bz\big]}(f)-\wz{[b, I_a^\bz\big]}(g)\r|
\le \lf|\wz{[b, I_a^\bz\big]}(f-g)\r|.$$
These combined with a standard density argument imply that
$\wz{[b, I_a^\bz]}$ admits a unique bounded extension from
$L^p(\gz)$ to $L^q(\gz)$.

Now we show (ii) by borrowing some ideas from  \cite{j78}. Since
$|z|^{n(1-\bz)}$ is infinitely differentiable in any open set away
from $0$, we choose $z_0\in\rn\setminus\{0\}$ and $\dz>0$ small
enough such that the function $|z|^{n(1-\bz)}$ can be expressed as
an absolutely convergent Fourier series in the neighborhood
$\{z\in\rn:\, |z-z_0|<2\dz\}$, that is,
$$|z|^{n(1-\bz)}=\sum_{m\in\zz^n} a_me^{i m\cdot z};$$
see, for example, \cite[p.\,266]{j78}.

Set $z_1\equiv \dz^{-1}z_0$.
Choose $\wz a\in(0, \fz)$ such that $a=\wz a(\wz a+1)(2+|z_1|)$.
For any $B\in\cb_{\wz a}$, we denote by $B'$
the ball centered at $c_B-r_Bz_1$ with radius $r_B$.
Notice that $$B'\subset B(c_B, (1+|z_1|)r_B)\subset B(c_{B'}, (1+2|z_1|)r_B).$$
This combined with \eqref{2.1} implies that $\gz(B)\sim\gz(B')$
with the equivalent constants depending on $z_1$.
Moreover, for any $x\in B$ and $y\in B'$, we have
$$\bigg|\dfrac{\dz(x-y)}{r_B}-z_0\bigg|
\le\bigg|\dfrac{\dz(x-c_B)}{r_B}\bigg|+\bigg|\dfrac{\dz(c_B-y)}{r_B}-z_0\bigg|
<2\dz.$$
Set $s(x)\equiv\mathrm{sgn} [b(x)-b_{B'}]$.
Notice that for any $y\in B'$ and $x\in B$, by \eqref{2.4}, we obtain
\begin{eqnarray*}
|y-x|&&\le|y-c_{B'}|+|c_{B'}-c_{B}|+|c_B-x|\\
&&<(2+|z_1|)r_B\le(2+|z_1|)\wz am(c_B)\le\wz a(\wz a+1)(2+|z_1|)m(x);
\end{eqnarray*}
thus, $B'\subset B(x, a m(x))$ for all $x\in B$.
Then we have
\begin{eqnarray*}
&&\dint_B|b(x)-b_{B'}|\,d\gz(x)\\
&&\hs=\dfrac{(r_B)^{n(1-\bz)}}{\gz(B')}\dint_B\dint_{B'}\dfrac{b(x)-b(y)}{[\dz|x-y|]^{n(1-\bz)}}
\lf|\dfrac{\dz(x-y)}{r_B}\r|^{n(1-\bz)}
s(x)\,d\gz(y)\,d\gz(x)\\
&&\hs=\dfrac{(\dz^{-1} r_B)^{n(1-\bz)}}{\gz(B')}\sum_{m\in\zz^n}a_m
\dint_\rn\dint_{B(x,\, am(x))}\dfrac{b(x)-b(y)}{[e^{-|x|^2}|x-y|^n]^{1-\bz}}
e^{i m\cdot \frac{\dz(x-y)}{r_B} }e^{-|x|^2(1-\bz)}\\
&&\hs\quad\times s(x)\chi_{B}(x)\chi_{B'}(y)\,d\gz(y)\,d\gz(x).
\end{eqnarray*}
If we set $f_m(y)\equiv e^{-i m\cdot \frac{\dz y}{r_B} }
\chi_{B'}(y)$ for all $y\in\rn$ and
$g_m(x)\equiv e^{i m\cdot \frac{\dz x}{r_B}}s(x)\chi_{B}(x)$ for all $x\in\rn$,
then the last formula above equals to a constant multiple of
\begin{eqnarray*}
\dfrac{(\dz^{-1}r_B)^{n(1-\bz)}}{\gz(B')}\sum_{m\in\zz^n}a_m
\dint_\rn\big[b, \wz{I}_a^\bz\big](f_m)(x)
e^{-|x|^2(1-\bz)}g_m(x)\,d\gz(x).
\end{eqnarray*}
Notice that $B'\subset B(c_B, (1+|z_1|)r_B)$ and \eqref{2.4} imply that
$e^{-|c_{B}|^2}\sim e^{-|c_{B'}|^2}$.
By this, $\supp g_m\subset B$ and the fact
$e^{-|x|^2}\sim e^{-|c_{B}|^2}$ for all $x\in B$,
we obtain that $e^{-|x|^2}\sim e^{-|c_{B'}|^2}$ for all $x\in B$.
From this and \eqref{2.2} together with H\"older's inequality,
we deduce that
\begin{eqnarray*}
\dint_B|b(x)-b_{B'}|\,d\gz(x)
&&\ls \dfrac{(\dz^{-1}r_B)^{n(1-\bz)}e^{-|c_{B'}|^2(1-\bz)}}
{\gz(B')}\sum_{m\in\zz^n}|a_m|
\dint_B \lf|[b, \wz{I}_{a}^\bz](f_m)(x)\r|\,d\gz(x)\\
&&\ls [\gz(B')]^{-\bz}[\gz(B)]^{1/q'}\sum_{m\in\zz^n}|a_m|
\lf\|[b, \wz{I}_{a}^\bz](f_m)\r\|_{L^q(\gz)}\\
&&\ls\lf\|[b, \wz{I}_{a}^\bz]\r\|_{L^p(\gz)\to L^q(\gz)} \gz(B).
\end{eqnarray*}
It follows that $\frac1{\gz(B)}\int_B|b(x)-b_B|\,d\gz(x)
\ls\|[b, \wz{I}_{a}^\bz]\|_{L^p(\gz)\to L^q(\gz)}$.
Taking the supermum over all balls $B\in\cb_a$ yields that
$\|b\|_\ast\ls \|[b, \wz{I}_{a}^\bz]\|_{L^p(\gz)\to L^q(\gz)}$.
Hence, we complete the proof of Theorem \ref{t1.1}.
\end{prf}

\newtheorem{proof}{\bf Proof of Theorem \ref{t1.2}}
\renewcommand\theproof{}

\begin{proof}\rm

To show (i),
by \eqref{2.2} and \eqref{2.3} together with an argument
similar to that used in the proof of Corollary \ref{c3.1}, we obtain
$$\lf|\wz{\big[b, \cm_a^\bz\big]}(f)(x)\r|
\ls\dint_{B(x,\, 2a(a+1)m(x))}
\dfrac{|b(x)-b(y)||f(y)|}{[e^{-|y|^2}|x-y|^n]^{1-\bz}}\,d\gz(y)
\sim T_{a^\ast}^\bz(b; f)(x),$$
where $a^\ast\equiv 2a(a+1)$.
Then following the argument used in the proof of Theorem \ref{t1.1} (i) yields that
$\wz{[b, \cm_a^\bz]}$ admits a unique bounded extension from
$L^p(\gz)$ to $L^q(\gz)$ with norm at most a constant multiple
of $\|b\|_\ast$.

Next we show (ii). Using H\"older's inequality
and the boundedness of $\wz{[b, \cm_a^\bz]}$, we obtain that for all $B\in\cb_a$,
\begin{eqnarray*}
\dfrac1{\gz(B)}\dint_B|b(y)-b_B|\,d\gz(y)
&&\le\dfrac1{[\gz(B)]^{1+\bz}}\dint_B
\dfrac1{[\gz(B)]^{1-\bz}}\dint_B|b(y)-b(x)|\chi_B(x)\,d\gz(x)\,d\gz(y)\\
&&\le\dfrac1{[\gz(B)]^{1+\bz}}\dint_B\wz{[b, \cm_a^\bz]}(\chi_B)(y)\,d\gz(y)\\
&&\le [\gz(B)]^{1/q'-1-\bz}\lf\|\wz{[b, \cm_a^\bz]}(\chi_B)\r\|_{L^q(\gz)}\le\lf\|\wz{[b, \cm_a^\bz]}\r\|_{L^p(\gz)\to L^q(\gz)},
\end{eqnarray*}
which together with the fact $b\in L^1(\gz)$ further implies that $b\in \BMO(\gz)$.
Thus, we complete the proof of Theorem \ref{t1.2}.
\end{proof}

\newtheorem{prof}{\bf Proof of Theorem \ref{t1.3}}
\renewcommand\theprof{}

\begin{prof}\rm
First we prove (i). Observe that for all $x\in\rn$,
\begin{equation}\label{4.10}
[b, \cm_a^\bz](f)(x)=b(x)\cm_a^\bz(f)(x)-\cm_a^\bz(bf)(x)\le\wz{[b, \cm_a^\bz]}(f)(x).
\end{equation}
Moreover, if we further assume that $b$ is a nonnegative function,
then for all $x\in\rn$,
\begin{equation}\label{4.11}
|[b, \cm_a^\bz](f)(x)|\le\wz{[b, \cm_a^\bz]}(f)(x).
\end{equation}
For all $x\in\rn$, since
\begin{eqnarray*}
-[b, \cm_a^\bz](f)(x)
&&\le \cm_a^\bz(b^+f)(x)+\cm_a^\bz(b^-f)(x)-b^+(x)\cm_a^\bz(f)(x)+b^-(x)\cm_a^\bz(f)(x)\\
&&\le |[b^+, \cm_a^\bz](f)(x)|+\cm_a^\bz(b^-f)(x)+b^-(x)\cm_a^\bz(f)(x),
\end{eqnarray*}
we then apply \eqref{4.11} and the fact
$\wz{[b^+, \cm_a^\bz]}(f)(x)\le \wz{[b, \cm_a^\bz]}(f)(x)$
to obtain
$$
 -[b, \cm_a^\bz](f)(x)
 \le \wz{[b, \cm_a^\bz]}(f)(x)+\cm_a^\bz(b^-f)(x)+b^-(x)\cm_a^\bz(f)(x),
$$
which together with \eqref{4.10} yields
$$
 |[b, \cm_a^\bz](f)(x)|\le \wz{[b, \cm_a^\bz]}(f)(x)
 +\cm_a^\bz(b^-f)(x)+b^-(x)\cm_a^\bz(f)(x).
$$
Then applying Theorem \ref{t1.2} (i) and Corollary \ref{c3.1}
together with the assumption $b^-\in L^\fz(\gz)$ yields
Theorem \ref{t1.3} (i).

To prove (ii), we first show $b\in\BMO(\gz)$.
To this end, fix $B\in\cb_a$. Then set
$$E\equiv \{x\in B:\, b(x)\le b_B\}.$$
Recall that $b_B$ is the integral average over the ball $B$.
It is not difficult to see that
$$\dint_{E}|b(x)-b_B|\,d\gz(x)
=\dint_{B\setminus E}|b(x)-b_B|\,d\gz(x).$$
Thus,
\begin{eqnarray}\label{4.12}
\dfrac{1}{\gz(B)}\dint_B|b(x)-b_B|\,d\gz(x)
&&=\dfrac{2}{\gz(B)}\dint_{E}[b_B-b(x)]\,d\gz(x)\\
&&=\dfrac{2}{[\gz(B)]^{2}}\dint_{E}\dint_B [b(y)-b(x)]
\,d\gz(y)\,d\gz(x).\noz
\end{eqnarray}
For all $x\in B$ and $B\in\cb_a$, we have
\begin{eqnarray*}
&&\sup_{\wz B\in\cb_a(x)}\dfrac{\gz(B\cap \wz B)}
{[\gz(\wz B)]^{1-\bz}[\gz(B)]^\bz}\\
&&\hs=\max\lf\{\sup_{\gfz{\wz B\in\cb_a(x)}{\gz(B)\ge\gz(\wz B)}}\dfrac{\gz(B\cap \wz B)}
{\gz(\wz B)}\lf(\dfrac{\gz(\wz B)}{\gz(B)}\r)^\bz,\,
\sup_{\gfz{\wz B\in\cb_a(x)}{\gz(B)<\gz(\wz B)}}\dfrac{\gz(B\cap \wz B)}
{\gz(B)}\lf(\dfrac{\gz( B)}{\gz(\wz B)}\r)^{1-\bz}\r\}
\le1,
\end{eqnarray*}
which further implies that
\begin{eqnarray*}
\sup_{\wz B\in\cb_a(x)}\dfrac{\gz(B\cap \wz B)}
{[\gz(\wz B)]^{1-\bz}[\gz(B)]^\bz}=1.
\end{eqnarray*}
Equivalently speaking, for all $B\in\cb_a$ and $x\in B$,
\begin{equation}\label{4.13}
\cm_a^\bz(\chi_B)(x)=[\gz(B)]^{\bz}.
\end{equation}
This combined with  \eqref{1.4} and \eqref{1.5} yields that for all $x\in \rn$,
\begin{eqnarray*}
\lf|\big[b, \cm_a^\bz\big](\chi_{B})(x)\r|
&&\ge \sup_{\wz B\in\cb_a(x)}\dfrac{1}{[\gz(\wz B)]^{1-\bz}}
\dint_{\wz B} b(y)\chi_{B}(y)\,d\gz(y)
-b(x)\cm_a^\bz(\chi_B)(x)\\
&&\ge \dfrac{1}{[\gz(B)]^{1-\bz}}
\dint_{B}[b(y)-b(x)]\,d\gz(y).
\end{eqnarray*}
Inserting this into \eqref{4.12} yields
\begin{eqnarray*}
\dfrac{1}{\gz(B)}\dint_B|b(x)-b_B|\,d\gz(x)
&&\le \dfrac{2}{[\gz(B)]^{1+\bz}}\dint_{E}
\lf|\big[b, \cm_a^\bz\big](\chi_{B})(x)\r|\,d\gz(x)\\
&&\le \dfrac{2}{[\gz(B)]^{1+\bz}}\dint_B
\lf|\big[b, \cm_a^\bz\big](\chi_{B})(x)\r|\,d\gz(x).
\end{eqnarray*}
Then using H\"older's inequality and the
boundedness of $[b, \cm_a^\bz]$ from $L^p(\gz)$ to $L^q(\gz)$, we further obtain that
for all $B\in\cb_a$,
\begin{eqnarray*}
\dfrac{1}{\gz(B)}\dint_B|b(x)-b_B|\,d\gz(x)
\le\dfrac{2[\gz(B)]^{1/q'}}{\gz(B)^{1+\bz}}\lf\|\big[b, \cm_a^\bz\big](\chi_{B})\r\|_{L^q(\gz)}\le2\lf\|\big[b, \cm_a^\bz\big]\r\|_{L^p(\gz)\to L^q(\gz)}.
\end{eqnarray*}
Taking the supremum over all balls $B\in\cb_a$ yields
$\|b\|_\ast^{\cb_a,\,1}\ls\|[b, \cm_a^\bz]\|_{L^p(\gz)\to L^q(\gz)}.$
This and the hypothesis $b\in
L^1(\gz)$ together with \eqref{2.8} and \eqref{1.1} imply that $b\in\BMO(\gz)$.

We still need to prove that $b^-\in L^\fz(\gz)$.
The differentiation theorem for the integral
implies that for almost all $x\in\rn$,
\begin{equation}\label{4.14}
b^+(x)=\lim_{B\in\cb_a(x),\, r_B\to0} \dfrac1{\gz(B)}\dint_B
b^+(z)\,d\gz(z).
\end{equation}
Fix $x$ satisfying \eqref{4.14}. Then for any given $\ez>0$, there
exists $B_0\in\cb_a(x)$ such that for all $B\in\cb_a(x)$ and
$B\subset B_0$, we have
$$\lf|b^+(x)-\frac1{\gz(B)}\int_B
b^+(z)\,d\gz(z)\r|<\ez,$$
and hence,
$$b^+(x)\le \frac1{\gz(B)}\int_B
b^+(z)\,d\gz(z)+\ez.$$
It follows that for $x$ satisfying
\eqref{4.14},
\begin{eqnarray}\label{4.15}
b^-(x)&&\le \frac1{\gz(B)}\int_B b^+(z)\,d\gz(z)+\ez-b^+(x)+b^-(x)\\
&&=\frac1{\gz(B)}\int_B b^+(z)\,d\gz(z)+\ez-b(x)
\le\frac1{\gz(B)}\int_B |b(z)|\,d\gz(z)-b(x)+\ez.\noz
\end{eqnarray}
For the above $B$ satisfying $B\in\cb_a(x)$ and $B\subset B_0$, we
use \eqref{4.13} to obtain that
\begin{eqnarray*}
\frac1{\gz(B)}\int_B |b(z)|\,d\gz(z)-b(x)
&&=\frac1{\gz(B)}\int_B |b(z)|\,d\gz(z)-b(x)\cm_a^\bz(\chi_B)(x)[\gz(B)]^{-\bz}\\
&&\le [\gz(B)]^{-\bz} \lf|[b, \cm_a^\bz](\chi_B)(x)\r|.
\end{eqnarray*}
Inserting this into \eqref{4.15} yields
\begin{equation}\label{4.16}
b^-(x)\le [\gz(B)]^{-\bz} \lf|[b, \cm_a^\bz](\chi_B)(x)\r|+\ez.
\end{equation}
We take the integration average over $B$ on both sides of \eqref{4.16},
then use H\"older's inequality and the boundedness of $[b,
\cm_a^\bz]$, and finally obtain
\begin{eqnarray*}
\dfrac1{\gz(B)}\dint_Bb^-(x)\,d\gz(x)
&&\le \dfrac1{\gz(B)^{1+\bz}} \dint_B\lf|[b, \cm_a^\bz](\chi_B)(x)\r|\,d\gz(x)+\ez\\
&&\le[\gz(B)]^{1/q'-1-\bz}\|[b, \cm_a^\bz](\chi_B)\|_{L^q(\gz)}+\ez
\le\|[b, \cm_a^\bz]\|_{L^p(\gz)\to L^q(\gz)}+\ez.
\end{eqnarray*}
This combined with the differentiation theorem for the integral implies that
for the points $y\in B_0$ such that
\eqref{4.14} holds with $b^+$ replaced by $b^-$,
$$b^-(y)\le \|[b, \cm_a^\bz]\|_{L^p(\gz)\to L^q(\gz)}+\ez.$$
So for the point $x$ satisfies \eqref{4.14} and also satisfies
\eqref{4.14} with $b^+$ replaced by $b^-$,
we have
$$b^-(x)\le \|[b, \cm_a^\bz]\|_{L^p(\gz)\to L^q(\gz)}+\ez.$$
Letting $\ez\to0$ yields
$$\|b^-\|_{L^\fz(\gz)}\le \|[b, \cm_a^\bz]\|_{L^p(\gz)\to L^q(\gz)}.$$
Thus, we obtain the desired results of (ii).
This finishes the proof of Theorem \ref{t1.3}.
\end{prof}

\bigskip

{\sc Liguang Liu} \& {\sc Dachun  Yang} (Corresponding author)

\medskip

School of Mathematical Sciences, Beijing Normal University,
Laboratory of Mathematics and Complex Systems, Ministry of
Education, Beijing 100875, People's Republic of China

\smallskip

{\it E-mails}: \texttt{liuliguang@mail.bnu.edu.cn} \&
\texttt{dcyang@bnu.edu.cn}

\end{document}